\documentclass[12pt]{article}
\usepackage{amssymb,amsmath,amsthm,enumerate}
\usepackage[dvips]{graphicx,xcolor}
\usepackage{here}

\newtheorem{theo}{Theorem}[section]
\newtheorem{coro}{Corollary}

\newtheorem{lem}[theo]{Lemma}

\newcommand{\pf}{\vspace{3mm}\noindent{\it Proof)}\hspace{5mm}}

\newcommand{\R}{{\bf R}}

\newcommand{\p}{\partial}

\newcommand{\lsm}{\lesssim}
\newcommand{\ep}{\varepsilon}

\newcommand{\al}{\alpha}
\newcommand{\bt}{\beta}

\newcommand{\lm}{\lambda}

\newcommand{\dsp}{\displaystyle}

\newcommand{\lr}{\longrightarrow }

\newcommand{\nb}{\nabla}
\newcommand{\fa}{\frac}
\newcommand{\sr}{\sqrt}
\newcommand{\lp}[1]{\left(#1 \right)}
\newcommand{\lb}[1]{\left\{#1 \right\}}
\newcommand{\lbt}[1]{\left[#1 \right]}

\begin{document}
 
\addtocounter{footnote}{1}
 \title{}
 
\begin{center}
{\large\bf 
On Glassey's conjecture for semilinear wave equations \\
in Friedmann-Lema\^itre-Robertson-Walker spacetime \\

}
\end{center}
\vspace{3mm}
\begin{center}

Kimitoshi Tsutaya$^\dagger$ 
        and Yuta Wakasugi$^\ddagger$ \\
\vspace{1cm}

 $^\dagger$Graduate School  of Science and Technology \\
Hirosaki University  \\
Hirosaki 036-8561, Japan\\
\footnotetext{* The research was supported by JSPS KAKENHI Grant Number JP18K03351. }

\vspace{5mm}
        $^\ddagger$
Graduate School of Engineering \\
Hiroshima University \\
Higashi-Hiroshima, 739-8527, Japan

\end{center}

\vspace{3mm}
\begin{abstract}
Consider nonlinear wave equations in the spatially flat Friedmann-Lema\^itre-Robertson-Walker (FLRW) spacetimes. We show blow-up in finite time of solutions and upper bounds of the lifespan of blow-up solutions to give the FLRW spacetime version of Glassey's conjecture for the time derivative nonlinearity. 
We also show blow-up results for the space time derivative nonlinearity. 
\end{abstract}  


\addtolength{\baselineskip}{2mm}

\vspace{1cm}
\section{Introduction.}

The spatially flat FLRW metric is given by 
\[
g: \; ds^2=-dt^2+a(t)^2d\sigma^2, 
\]
where the speed of light is equal to $1$, $d\sigma^2$ is the line element of $n$-dimensional Euclidean space and $a(t)$ is the scale factor, which describes expansion or contraction of the spatial metric. 
As in our earlier work \cite{TW1,TW2,TW3}, we treat the scale factor as 
\begin{equation}
a(t)=ct^{\fa 2{n(1+w)}}  
\label{scw} 
\end{equation}
where $c$ is a positive constant, and $w$ is the proportionality constant in the range $-1< w\le 1$. The constant $w$ appears in the equation of state relating the pressure to the density for the perfect fluid. See \cite{TW1}. 

In the preceding papers \cite{TW1,TW2,TW3}, we have shown upper bounds of the lifespan
for the equation $\Box_g u=-|u|^p$. 
In this paper we consider the equation $\Box_g u=-|u_t|^p$ as well as $\Box_g u=-|\nb_x u|^p$, where $\nb_xu = (\p u/\p x^1, \cdots,\p u/\p x^n), \; (x^1,\cdots, x^n)\in \R^n$. 

For the flat FLRW metric with \eqref{scw},  
the semilinear wave equation $\Box_g u\\
=|g|^{-1/2}\p_\al(|g|^{1/2}g^{\al\bt}\p_\bt)u= -|u_t|^p$ or $-|\nb_x u|^p$ with $p>1$ becomes 
\begin{equation}
u_{tt}-\fa 1{t^{4/(n(1+w))}}\Delta u+\fa 2{(1+w)t}u_t=|u_t|^p, \mbox{ or }|\nb_x u|^p,   \quad  x\in \R^n  
\label{ore}
\end{equation}
where  $\Delta=\p_1^2+\cdots \p_n^2, \; \p_j=\p/\p x^j, \; j=1,\cdots,n$.   
Our aim of this paper is to show that blow-up in a finite time occurs for the equation above as well as upper bounds of the lifespan of the blow-up solutions.

We first consider the following Cauchy problem
in order to compare with the related known results including the case of the Minkowski spacetime: 
\begin{equation}
u_{tt}-\fa 1{t^{2\al}}\Delta u+\fa\mu{t}u_t=|u_t|^p, \qquad  t>1, \; x\in \R^n 
\label{Prob0}
\end{equation}
with the initial data given at $t=1$, 
\begin{equation}
u(1,x)=\ep u_0(x), \; u_t(1,x)=\ep u_1(x), \qquad x\in \R^n, 
\label{data0}
\end{equation}
where $\al$ and $\mu$ are nonnegative constants and $\ep>0$ is a small parameter.  

Let $T_\ep$ be the lifespan of solutions of \eqref{Prob0} and \eqref{data0}, say, $T_\ep$ is the supremum of $T$ such that  \eqref{Prob0} and \eqref{data0} have a solution for $x\in \R^n$ and $1\le t<T$. 

Let $\al=\mu=0$ and $p_G(n)=1+2/(n-1)$.  
The so-called Glassey's conjecture \cite{Gl2} asserts that if $p>p_G(n)$, then there exist global solutions in time for small initial data, on the other hand, if $1<p\le p_G(n)$ with $n\ge 2$ or if $p>1$ for $n=1$, then blow-up in finite time occurs. 
This conjecture is proved to be almost true. Actually, blow-up results in low dimensions ($n=2,3$), or in high dimensions ($n\ge 4$) imposing radial symmetry were proved in, e.g., \cite{Ag,Gl1,Jo2,Ma,Ra,Sch}, and Zhou \cite{Zh} finally gave a simple proof of the blow-up result for $1<p\le p_G(n)$ and $n\ge 2$ as well as for $p > 1$ and $n=1$. 
Global existence of solutions in low dimensions ($n=2,3$) has been proved in, e.g., \cite{HT,Si,Tz}. 
For high dimensions ($n\ge 3$), it is proved by \cite{HWY} that there exist global solutions in the radial case for $p>p_G(n)$. They \cite{HWY} also proved the lifespan of local solutions in time for $1<p\le p_G(n)$.

For the case $\al=0$ and $\mu\ge 0$, 
it is recently shown by Hamouda and Hamza \cite{HH}  
that blow-up in finite time occurs  and the lifespan of the blow-up solutions satisfies 
\begin{align}
&T_\ep \le C\ep^{-(p-1)/\{1-(n+\mu-1)(p-1)/2\}}&& \mbox{ if }\quad  1< p < p_G(n+\mu), \; n\ge 1, 
\label{subst}\\
&T_\ep \le \exp(C\ep^{-(p-1)})&&\mbox{ if }\quad p = p_G(n+\mu), \; n\ge 1. 
\label{crst}
\end{align}
These results improve the ones in \cite{LT}. 

The present paper treats the case $\al\ge 0$ and $\mu\ge 0$. 
We first show blow-up in a finite time and upper estimates of the lifespan of solutions of \eqref{Prob0} and \eqref{data0} in the case $0\le \al<1$. 
If $\al=0$, our upper bounds of the lifespan coincide with the results above by \cite{HH}. 
Similar results are independently shown by \cite{HHP} where energy solutions are treated.
In our results, however, another exponent appears as a blow-up condition in some case. 
This is different from the results  by \cite{HHP}. 
We emphasize that the generalized exponent of $p_G(n+\mu)$ cannot always be the critical exponent for the global existence of solutions. 
Our proofs are based on the test function method with 
the modified Bessel function of the second kind and on a generalized Kato's lemma.  
We next treat the case $\al\ge 1$. 

Moreover, we show blow-up results for the problem
\begin{equation}
\begin{cases}
\dsp u_{tt}-\fa 1{t^{2\al}}\Delta u+\fa\mu{t}u_t=|\nb_xu|^p, \qquad  t>1, \; x\in \R^n, \\
u(1,x)=\ep u_0(x), \; u_t(1,x)=\ep u_1(x), \qquad x\in \R^n. 
\end{cases}
\label{Prob0x}
\end{equation} 
Unlike the above equation \eqref{Prob0}, our blow-up conditions are related to exponents that originate from the Strauss and Fujita ones and to another exponent like in the case of the time derivative nonlinearity. Hence, upper bounds of the lifespan have to do with those exponents. 

We then apply our results for \eqref{Prob0} and \eqref{Prob0x} to the original equation \eqref{ore}. 
Our aim is especially to clarify the difference with the case of the Minkowski spacetime and also how the scale factor affects the lifespan of the solution.  
Since global existence of solutions has not been obtained yet, {\it critical} exponent used in this paper means a candidate of the true critical exponent.

The paper is organized as follows. 
In Section 2, we state our first main result for \eqref{Prob0} and \eqref{data0} in the case  $0\le \al <1$. 
Theorem \ref{th21} presents that it is possible in some case to improve the estimate of the lifespan affected by Glassey's exponent. 
To prove the theorem, we use the test function method, and also a generalized Kato's lemma for a first-order differential inequality, which is applied to the wave equation with the scale-invariant damping. This is proved by John's iteration argument \cite{Jo1}. 
We then show our second result which is for the case $\al\ge 1$. 
In Section 3 we treat \eqref{Prob0x} and divide results into several cases $0\le \al<1$, wavelike and heatlike cases, critical and subcritical cases, and $\al\ge 1$. 
Finally in Section 4, we apply the theorems in Sections 2 and 3 to the original equation \eqref{ore}. We discuss the effect of the scale factor to the solutions.

\vspace{5mm}

\section{Time derivative nonlinearity.}

\subsection{Case $0\le \al <1$.}

\setcounter{equation}{0}

We first consider the problem \eqref{Prob0} with \eqref{data0} for $0\le \al <1$. 
Our first result is the following theorem: 

\begin{theo} 
Let $n\ge 2, \; 0\le \al<1, \; \mu\ge 0$ and 
\begin{align*}
&1<p\le p_G'(n,\al,\mu)\equiv 1+\fa 2{(1-\al)(n-1)+\mu+\al},  
\intertext{ or }  
&1<p<p_0(n,\al,\mu)\equiv 1+\fa 1{n(1-\al)+\mu}.
\end{align*}
Assume that $u_0\in C^2(\R^n)$ and $u_1\in C^1(\R^n)$ are nontrivial and satisfy $u_1(x)\ge u_0(x)\ge 0$, $\mbox{\rm supp }u_0, \mbox{\rm supp }u_1\subset \{|x|\le R\}$ with $R>0$. 
Suppose that the problem \eqref{Prob0} with \eqref{data0} has  a classical solution $u\in C^2([1,T)\times\R^n)$.  
Then, $T<\infty$ and there exists a constant $\ep_0>0$ depending on 
$p,\al,\mu,R,u_0,u_1$ such that $T_\ep$ has to satisfy 
\begin{align}
T_\ep&\le C\ep^{\fa{-(p-1)}{1-\{(1-\al)(n-1)+\mu+\al\}(p-1)/2}} && \quad \mbox{if }
p<p_G'(n,\al,\mu), 
\label{subcls-21} \\
T_\ep&\le \exp(C\ep^{-(p-1)}) &&\quad \mbox{if }
p=p_G'(n,\al,\mu), 
\label{cls-21} \\
T_\ep&\le C\ep^{\fa{-(p-1)}{1-(p-1)\{n(1-\al)+\mu\}}} && \quad \mbox{if }
p<p_0(n,\al,\mu) 
\label{subcls-22}
\end{align}
for $0<\ep\le \ep_0$, 
where $C>0$ is a constant independent of $\ep$. 
\label{th21}
\end{theo}

\noindent
{\bf Remark}
 (1) If $\al =0$, then \eqref{subcls-21} and \eqref{cls-21} are the same with the upper bounds \eqref{subst} and \eqref{crst}. \\
(2) By the theorem, the exponent $p_G'(n,\al,\mu)$ cannot always be the critical exponent for the global existence of solutions. We discuss more details in the end of this subsection. \\
(3) If $p<p_0(n,\al,\mu)$, then the above assumption $u_1(x)\ge u_0(x)\ge 0$ can be replaced just by $u_1(x)\ge 0$. 

\pf
Mutiplying \eqref{Prob0} by a test funtion $\phi(t,x)$ and $t^\mu$, and 
integrating over $\R^n$, we have 
\begin{equation}
\fa d{dt}\int t^\mu(u\phi)_tdx-2\fa d{dt}\int t^\mu u\phi_t dx
+\int t^\mu u\lp{\phi_{tt}-\fa 1{t^{2\al}}\Delta\phi+\fa\mu t\phi_t}dx
=\int t^\mu |u_t|^p\phi dx. 
\label{9}
\end{equation}
Integrating over $[1,t]$, we obtain 
\begin{align}
&\int t^\mu(u\phi)_t(t,x)dx-2\int t^\mu u\phi_t(t,x) dx 
+\int_1^t s^\mu \int  u\lp{\phi_{ss}-\fa 1{s^{2\al}}\Delta\phi+\fa\mu s\phi_s}dx ds\nonumber\\
&=\ep\int (u_1(x)\phi(1,x)- u_0(x)\phi_t(1,x)) dx+\int_1^t s^\mu\int  |u_s|^p\phi dxds.
\label{F10}
\end{align}
We remark that 
the $C^2$-solution $u$ of \eqref{Prob0} and \eqref{data0} has the property of finite speed of propagation, and satisfies
\begin{equation}
\mbox{supp }u(t,\cdot)\subset \{|x|\le A(t)+R\},  
\qquad A(t)=\int_1^t s^{-\al}ds =\fa{t^{1-\al}-1}{1-\al}, 
\label{suppu0}
\end{equation}
provided that supp $u_0$, supp $u_1\subset \{|x|\le R\}$. See \cite{TW1} for its proof. 

We now define a smooth test function by
\begin{align*}
&\phi(t,x)=\lm(t)\int_{|\omega|=1}e^{x\cdot \omega}dS_\omega,  \\
&\lm(t)=t^{(1-\mu)/2}K_\nu\lp{\fa 1{1-\al}t^{1-\al}}, \qquad \nu=\fa{\mu-1}{2(1-\al)}, 
\end{align*}
where $K_\nu(t)$ is the modified Bessel funtion of the second kind which is given by 
\begin{equation}
K_\nu(t)=\int_0^\infty e^{-t\cosh z}\cosh \nu zdz, \qquad t>0, \quad \nu\in\R. 
\label{Bessel}
\end{equation}
It is well-known that the Bessel function $K_\nu$ satisfies the following properties (see, e.g.,\cite{AS}): 
\begin{align}
&t^2K_\nu''(t)+tK_\nu'(t)-(t^2+\nu^2)K_\nu(t)=0,   
  \label{B0}\\
& K_\nu(t)=\sr{\fa\pi 2}\fa{e^{-t}}{\sr t}\lp{1+O\lp{\fa 1t}} \qquad (t\to\infty),  
    \label{B1}\\
& K_\nu'(t)=\fa\nu tK_\nu(t)-K_{\nu+1}(t).  
   \label{B2}
\end{align}  
We can verify by \eqref{B0}-\eqref{B2} that there holds  
\begin{align}
&\phi_{tt}-\fa 1{t^{2\al}}\Delta\phi+\fa\mu t\phi_t=0. 
 \label{Ph1}
\end{align}  
The following estimate is shown in \cite{TW1}:
\begin{equation}
\int_{|x|\le A(t)+R}|\phi(t,x)|dx\lsm (t+R)^{\{(1-\al)(n-1)-(\mu-\al)\}/2}. 
\label{V}
\end{equation}
See (3.23) in \cite{TW1}. 
We also see by \cite{TW1} that 
\begin{equation}
C_d\equiv \ep\int (u_1(x)\phi(1,x)- u_0(x)\phi_t(1,x)) dx
=C_0\ep >0
\label{Cdata}
\end{equation}
under assumption on the initial data. 
Moreover, we have 
\begin{align}
\phi_t(t,x)&=\lm'(t)\int_{|\omega|=1}e^{x\cdot\omega}dS_\omega \nonumber\\
&=\lb{\fa{1-\mu}2 t^{-1}+t^{-\al}\fa{K_\nu'\lp{\fa 1{1-\al} t^{1-\al}}}{K_\nu\lp{\fa 1{1-\al} t^{1-\al}}}}\phi(t,x)\\
&=-t^{-\al}\fa{K_{\nu+1}\lp{\fa 1{1-\al} t^{1-\al}}}{K_\nu\lp{\fa 1{1-\al} t^{1-\al}}}\phi(t,x). 
\label{phit}
\end{align}
where we have used \eqref{B2} and $\nu=(\mu-1)/(2(1-\al))$ for the last equality. 

Set 
\[
F_1(t)=\int u\phi(t,x)dx.
\] 
From \eqref{F10}, proceeding as in \cite{TW1}, we have
\begin{align*}
\fa{t^{\mu-1}}{K_\nu\lp{\fa 1{1-\al}t^{1-\al}}^2}F_1(t)
-\fa 1{K_\nu\lp{\fa 1{1-\al}}^2}F_1(1)
&\ge C_d\int_1^t \fa{s^{-1}}{K_\nu\lp{\fa 1{1-\al}s^{1-\al}}^2}ds. 
\end{align*}
Since $F_1(1)/{K_\nu\lp{\fa 1{1-\al}}^2}>0$ by assumption, we obtain
\begin{equation}
F_1(t)=\int u\phi(t,x)dx>0 \quad \mbox{for }t\ge 1. 
\label{F_1+}
\end{equation}
We next go back to \eqref{9}, which becomes 
\[
\fa d{dt}\int t^\mu u_t\phi dx-\int t^\mu u_t\phi_t dx
-\fa 1{t^{2\al}}\int t^\mu u\Delta\phi dx
=\int t^\mu |u_t|^p\phi dx. 
\]
Using \eqref{phit} and $\Delta\phi=\phi$ yields
\begin{equation}
\fa d{dt}\int t^\mu u_t\phi dx+t^{-\al}\fa{K_{\nu+1}\lp{\fa 1{1-\al} t^{1-\al}}}{K_\nu\lp{\fa 1{1-\al} t^{1-\al}}}\int t^\mu u_t\phi dx
-\fa 1{t^{2\al}}\int t^\mu u\phi dx
=\int t^\mu |u_t|^p\phi dx. 
\label{91}
\end{equation}
On the other hand, by \eqref{F10}, \eqref{Ph1}, \eqref{Cdata} and \eqref{phit}, 
\begin{equation}
\int t^\mu u_t\phi(t,x)dx+t^{-\al}\fa{K_{\nu+1}\lp{\fa 1{1-\al} t^{1-\al}}}{K_\nu\lp{\fa 1{1-\al} t^{1-\al}}}\int t^\mu u\phi(t,x) dx 
=C_d+\int_1^t s^\mu\int  |u_s|^p\phi dxds.
\label{F}
\end{equation}
We see that there exists a constant $M\ge 1$ such that 
\begin{equation}
\fa{K_{\nu+1}\lp{\fa 1{1-\al} t^{1-\al}}}{K_\nu\lp{\fa 1{1-\al} t^{1-\al}}}\le M
\label{keyc}
\end{equation}
for $t\ge 1$ by \eqref{B1}. 
Combining \eqref{91} and \eqref{F} multiplied by $M^{-1}t^{-\al}$, we have
\begin{align*}
&\fa d{dt}\int t^\mu u_t\phi dx+\lbt{\fa{K_{\nu+1}\lp{\fa 1{1-\al} t^{1-\al}}}{K_\nu\lp{\fa 1{1-\al} t^{1-\al}}}+M^{-1}}t^{-\al}\int t^\mu u_t\phi dx\\
&\qquad +
\lbt{M^{-1}\fa{K_{\nu+1}\lp{\fa 1{1-\al} t^{1-\al}}}{K_\nu\lp{\fa 1{1-\al} t^{1-\al}}}-1}t^{-2\al}\int t^\mu u\phi(t,x) dx \\
=&C_dM^{-1}t^{-\al}+ \int t^\mu |u_t|^p\phi dx+M^{-1}t^{-\al}\int_1^t s^\mu\int  |u_s|^p\phi dxds.
\end{align*}
We note that 
\[
t^{-2\al+\mu}\int u\phi(t,x) dx >0
\]
for all $t\ge 1$ by \eqref{F_1+}. 
Using \eqref{keyc}, we obtain 
\begin{align}
&\fa d{dt}\int t^\mu u_t\phi dx+\lp{M+\fa 1M}t^{-\al}\int t^\mu u_t\phi dx \nonumber\\
\ge &C_dM^{-1}t^{-\al}+ \int t^\mu |u_t|^p\phi dx+M^{-1}t^{-\al}\int_1^t s^\mu\int  |u_s|^p\phi dxds
\qquad \mbox{for }t\ge 1.  
\label{star}
\end{align}
We now set
\[
G(t)=\int t^\mu u_t\phi dx-\fa 1{M^2+1}\int_1^t s^\mu\int  |u_s|^p\phi dxds-\fa{C_d}{2(M^2+1)} \qquad \mbox{for }t\ge 1. 
\]
Then 
\[
G'(t)=\fa{d}{dt}\int t^\mu u_t\phi dx-\fa 1{M^2+1} \int t^\mu |u_t|^p\phi dx, 
\]
hence \eqref{star} becomes 
\[
G'(t)+\fa{M^2+1}{M}t^{-\al}G(t)\ge \fa{C_d}{2M}t^{-\al}+\fa{M^2}{M^2+1}\int t^\mu |u_t|^p\phi dx>0 \quad \mbox{for }t\ge 1. 
\]
Multiplying the inequality above by $\exp[(M^2+1)t^{1-\al}/(M(1-\al))]$ and integrating over $[1,t]$, we obtain
\begin{equation}
\exp\lp{\fa{M^2+1}{M(1-\al)}t^{1-\al}}G(t)-\exp\lp{\fa{M^2+1}{M(1-\al)}}G(1)>0
\quad \mbox{for }t\ge 1. 
\label{G1}
\end{equation}
Note that 
\begin{align*}
G(1)&=\ep\int u_1(x)\phi(1,x)dx-\fa \ep{2(M^2+1)}\int (u_1(x)\phi(1,x)- u_0(x)\phi_t(1,x)) dx\\
&\ge \ep\int\lp{\fa{2M^2+1}{2(M^2+1)}u_1(x)-\fa M{2(M^2+1)}u_0(x)}\phi(1,x)dx
>0 
\end{align*}
by \eqref{phit}, \eqref{keyc} and assumption on the initial data. It holds from \eqref{G1} that $G(t)>0$ for $t\ge 1$. Thus, we see that 
\begin{equation}
\int t^\mu u_t\phi dx>\fa 1{M^2+1}\int_1^t s^\mu\int  |u_s|^p\phi dxds+\fa{C_d}{2(M^2+1)} \qquad \mbox{for }t\ge 1. 
\label{star1}
\end{equation}
We now define 
\[
H(t)=\int_1^t s^\mu\int  |u_s|^p\phi dxds+\fa{C_d}2
\quad \mbox{for }t\ge 1. 
\]
By H\"older's inequality, \eqref{V} and \eqref{star1}, we have 
\begin{align}
H'(t)=\int t^\mu |u_t|^p\phi dx
&\ge \lp{\int t^\mu u_t\phi dx}^p\lp{\int_{|x|\le A(t)+R} t^\mu \phi dx}^{1-p} \nonumber\\
&\gtrsim t^{-\{(1-\al)(n-1)+\mu+\al\}(p-1)/2}H(t)^p 
\quad \mbox{for }t\ge 1. 
\label{Ht}
\end{align}
We also have 
\begin{equation}
H(1)=\fa{C_d}2=\fa 12C_0\ep>0
\label{H1}
\end{equation}
by \eqref{Cdata}. 
By integrating \eqref{Ht} multiplied by $H(t)^{-p}$ over $[1,t]$ and using \eqref{H1}, we therefore obtain the desired results \eqref{subcls-21} and \eqref{cls-21}. 


It remains to prove \eqref{subcls-22}. We use the following lemma, which is a generalized Kato's lemma. 

\begin{lem}
Let $p>1, \;a\ge 0, \;b\ge 0, \;q\ge 0, \; r\ge 0, \; \mu\ge 0, \; c>0$ and 
\[
M\equiv (p-1)(c-a)-q+1>0. 
\]
Let $T\ge T_1>T_0\ge 1$. 
Assume that $F\in C^1([T_0,T))$ satisfies the following three conditions:  
\begin{align*}
(i) \quad & F(t) \ge A_0t^{-a}(\ln t)^{-b}(t-T_1)^c \qquad \mbox{for }t\ge T_1, \\
(ii) \quad &F'(t) +\fa{\mu}{t}F(t)\ge A_1(t+R)^{-q}(\ln t)^{-r}|F(t)|^p \quad \mbox{for }t\ge T_0,	 \\
(iii) \quad &F(T_0)> 0,  
\end{align*}
where $A_0, A_1$ and $R$ are positive constants. Then,  
$T$ has to satisfy
\[
T^{M/(p-1)}(\ln t)^{-b-r/(p-1)}<CA_0^{-1},  
\]
where $C$ is a constant depending on $R,A_1,\mu,p,q,r,a,b$ and $c$. 
\label{lm22}
\end{lem}

\pf 
Mutiplying assumption (ii) by $t^\mu$, we have  
\[
t^\mu F'+\mu t^{\mu-1}F\ge A_1 t^{\mu}(t+R)^{-q}(\ln t)^{-r}|F|^p. 
\]
Integrating the above inequality over $[T_0,t]$ yields 
\begin{align}
t^\mu F(t)-F(T_0)\ge A_1C_{R,q}\int_{T_0}^t s^{\mu-q}(\ln s)^{-r}|F(s)|^pds\ge 0, \quad C_{R,q}=(1+R)^{-q}. 
\label{e0}
\end{align}
By assumption (iii),  we see that $F(t)>0$ for $t\ge T_0$. 
Hence, by assumption (i), we have
\begin{align*}
F(t) &\ge A_0^pA_1C_{R,q} t^{-\mu}\int_{T_1}^t s^{\mu-ap-q}(\ln s)^{-bp-r}(s-T_1)^{cp} ds\\
 &\ge A_0^pA_1C_{R,q} t^{-\mu-ap-q}(\ln t)^{-bp-r}  \int_{T_1}^t (s-T_1)^{\mu+cp}ds \\
&= \fa{A_0^pA_1C_{R,q}}{\mu+cp+1}t^{-\mu-ap-q}(\ln t)^{-bp-r}(t-T_1)^{\mu+cp+1}
\qquad \mbox{for }t\ge T_1. 
\end{align*} 
Based on the fact above, we define the sequences $a_j, \; b_j, \; c_j, \; D_j$ for $j = 0, 1, 2, \cdots$ by
\begin{align}
a_{j+1} = pa_j + \mu+q, & & b_{j+1} = pb_j+r,  
& & c_{j+1}=pc_j +\mu+1, & & D_{j+1} = \frac{A_1C_{R,q}D_j^p}{pc_j+\mu+1}  \label{e2}\\
a_0 = a, & & b_0 = b, & & c_0 = c, & & D_0= A_0.  \label{e3}
\end{align}
Solving \eqref{e2} and \eqref{e3}, we obtain  
\begin{align*}
&a_j = p^j\lp{a+\fa{\mu+q}{p-1}}-\fa{\mu+q}{p-1}, \qquad 
b_j = p^j\lp{b+\fa{r}{p-1}}-\fa{r}{p-1},  \\
&c_j = p^j\lp{c+\fa{\mu+1}{p-1}}-\fa{\mu+1}{p-1}, 
\end{align*}
and thus
\[
D_{j+1} = \frac{A_1C_{R,q}D_j^p}{c_{j+1}} \ge \lp{c+\fa{\mu+1}{p-1}}^{-1}\frac{A_1C_{R,q}D_j^p}{p^{j+1}}.
\]
Then, 
\begin{align*}
D_j &\ge \frac{BD_{j-1}^p}{p^{j}} \\
&\ge \fa B{p^{j}}\lp{\fa{BD_{j-2}^p}{p^{j-1}}}^p =\fa{B^{1+p}}{p^{j+p(j-1)}}D_{j-2}^{p^2} \\
&\ge  \fa{B^{1+p}}{p^{j+p(j-1)}}\lp{\fa{BD_{j-3}^p}{p^{j-2}}}^{p^2}
=\fa{B^{1+p+p^2}}{p^{j+p(j-1)+p^2(j-2)}}D_{j-3}^{p^3} \\
&\ge \cdots\cdots \ge \fa{B^{1+p+p^2+\cdots+p^{j-1}}}{p^{j+p(j-1)+p^2(j-2)+\cdots+p^{j-1}}}D_0^{p^j}, \\
\intertext{and }
\ln D_j &\ge 
 \frac{\ln B}{p-1}(p^j-1)-p^j \sum_{k=0}^j 
\frac{k}{p^k}\ln p+p^j\ln D_0, 
\end{align*}
where $B=\{c+(\mu+1)/(p-1)\}^{-1}A_1C_{R,q}$.  
For sufficiently large $j$, we have
\[
D_j \ge \exp(Ep^j),
\]
where
\begin{equation}
E = \frac{1}{p-1}\min\left(0, \ln B \right) - \sum_{k=0}^{\infty}\frac{k}{p^k}\ln p+\ln A_0. 
\label{eqofE}
\end{equation}
Thus, since $F(t)\ge D_jt^{-a_j}(\ln t)^{-b_j}(t-T_1)^{c_j}$ holds for $t\ge T_1$, we obtain
\begin{align}
F(t)  &\ge t^{(\mu+q)/(p-1)}(\ln t)^{r/(p-1)}(t-T_1)^{-(\mu+1)/(p-1)} \nonumber\\
  & \quad \cdot\exp\left[\left\{E+\lp{c+\fa{\mu+1}{p-1}}\ln(t-T_1)-\lp{a+\fa{\mu+q}{p-1}}\ln t-\lp{b+\fa{r}{p-1}}\ln(\ln t)\right\}p^j\right]
 \label{5}
\end{align}
for $t\ge T_1$. 
Since     
\[
\lp{c+\fa{\mu+1}{p-1}}-\lp{a+\fa{\mu+q}{p-1}}=c-a+\fa{1-q}{p-1}>0
\]
by assumption, choosing $t$ large enough, we can find a positive $\delta$ such that
\[
E+\lp{c+\fa{\mu+1}{p-1}}\ln(t-T_1)-\lp{a+\fa{\mu+q}{p-1}}\ln t-\lp{b+\fa{r}{p-1}}\ln(\ln t) \ge \delta > 0. 
\]
It then follows from \eqref{5}  that $F(t) \lr \infty$ as $j \to \infty$ for 
sufficiently large $t$. We therefore see that the lifespan $T$ of $F(t)$ has to satisfy
\[
T^{M/(p-1)}(\ln T)^{-b-r/(p-1)}<CA_0^{-1},  
\]
where $M=(p-1)(c-a)-q+1>0$, and $C$ is a constant depending on $A_1,R,\mu,p,q,r,a,b$ and $c$. 
This completes the proof of the proposition. 
\hfill\qed 

We now prove \eqref{subcls-22} by applying Lemma \ref{lm22}. 
Set  
\begin{equation}
F(t)=\int u_t(t,x)dx.
\label{int-u_t}
\end{equation}
Integrating the equation \eqref{Prob0} and using H\"older's inequality, we have by \eqref{suppu0}, 
\begin{align}
F'(t)+\fa\mu t F(t) &=\int |u_t|^p dx  
   \label{eqF} \\
 & \ge \fa 1{(A(t)+R)^{n(p-1)}}|F(t)|^p \nonumber \\
 & \ge \fa C{(t+R)^{n(1-\al)(p-1)}}|F(t)|^p. 
\label{Ho1}
\end{align}
Mutiplying \eqref{Ho1} by $t^\mu$
 and integrating imply
\begin{equation}
t^\mu F(t)-F(1)\gtrsim \int_1^t s^{\mu-n(1-\al)(p-1)} |F(s)|^pds\ge 0. 
\label{F01}
\end{equation}
We hence see that 
\begin{equation} 
F(t)\ge t^{-\mu}F(1)=C\ep t^{-\mu}>0 \qquad \mbox{for }t\ge 1
\label{F2}
\end{equation}
by assumption. 
From \eqref{F01} and \eqref{F2}, 
we have 
\begin{align}
F(t)&\ge C\ep^p t^{-\mu}\int_1^t s^{\mu-n(1-\al)(p-1)-\mu p}ds  \nonumber \\
&\ge C\ep^p t^{-\mu(p+1)-n(1-\al)(p-1)}\int_1^t (s-1)^\mu ds.  \nonumber\\
\intertext{Therefore, we obtain}
F(t)&\ge C\ep^p t^{-\mu(p+1)-n(1-\al)(p-1)}(t-1)^{\mu+1} 
\qquad \mbox{for }t\ge 1. 
\label{ann}
\end{align}
Finally, by \eqref{Ho1} and \eqref{ann}, applying Lemma \ref{lm22} with $q=n(1-\al)(p-1)$, $a=\mu(p+1)+n(1-\al)(p-1)$, $b=r=0, \; c=\mu+1$ and $A_0=C\ep^p$ 
, we obtain the desired result \eqref{subcls-22} since
\begin{align*}
M &=(p-1)\lb{\mu+1 - \mu(p+1)-n(1-\al)(p-1)}-n(1-\al)(p-1)+1  \\
&=p\lb{1-(p-1)(\mu+n(1-\al))}>0. 
\end{align*}
This completes the proof of Theorem \ref{th21}. 
\hfill\qed

\begin{figure}[h!]
\includegraphics[width=13cm, bb=120 570 500 760, clip]{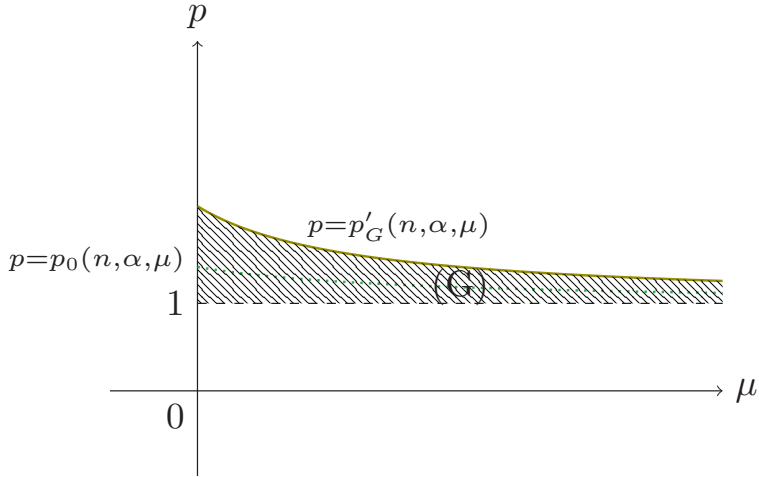}
\caption{\label{Fig1}Range of blow-up conditions in case $n=3$ and $\alpha=0.2$}
\end{figure}

\begin{figure}[h!]
\includegraphics[width=15cm, bb=120 550 500 760, clip]{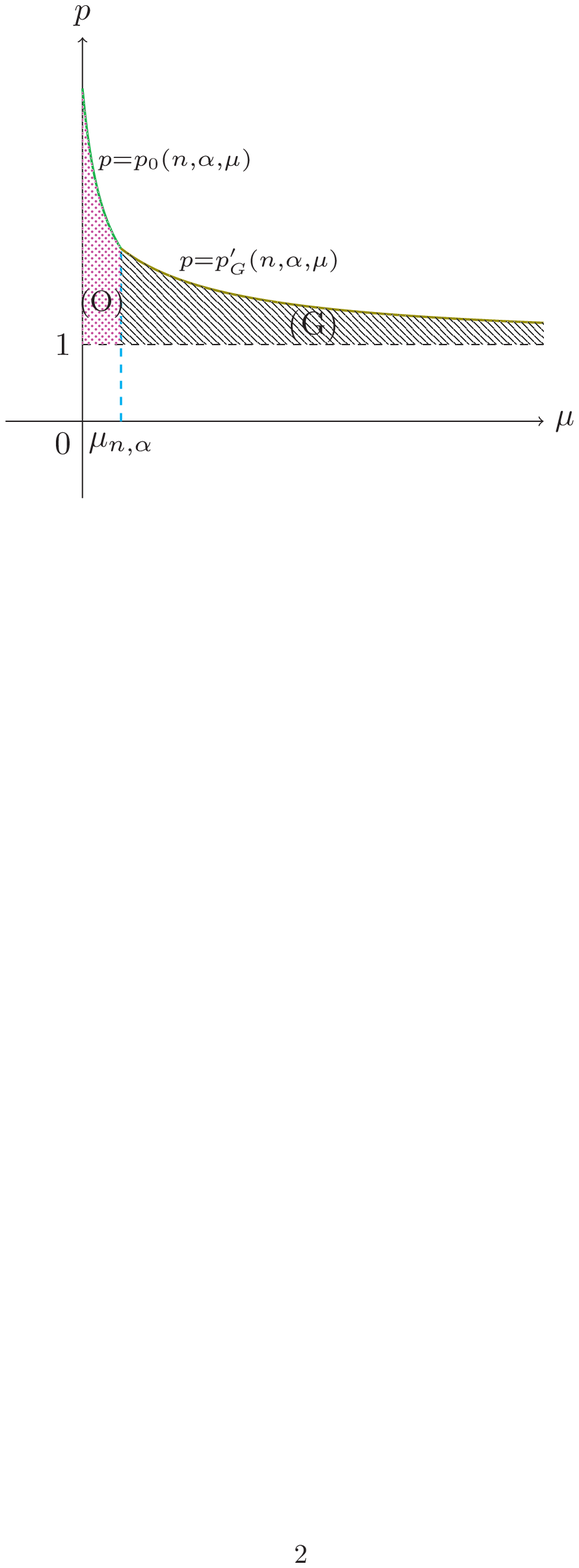}
\caption{\label{Fig2}Range of blow-up conditions in case $n=3$ and $\alpha=0.9$}
\end{figure}

In the end of this subsection, we discuss the blow-up condtions and the estimates of the lifespan in the two subcritical cases in Theorem \ref{th21}. 
We note that if $p_G'(n,\al,\mu)=p_0(n,\al,\mu)$, then 
\[
\mu=\mu_{n,\al}\equiv \al(n+2)-(n+1). 
\]
Fig.~\ref{Fig1}
and Fig.~\ref{Fig2} 
below show the regions of blow-up conditions in the cases  $n=3, \; \al=0.2$ and $n=3, \; \al=0.9$, respectively.

For $\mu > \max\{0,\mu_{n,\al}\}$, the exponent $p_G'(n,\al,\mu)$ is bigger than 
$p_0(n,\al,\mu)$. See the region (G) in Fig.s~\ref{Fig1} and \ref{Fig2}. 
On the other hand, if $(n+1)/(n+2)\le \al<1$, then for 
$0\le \mu \le \mu_{n,\al}$, the condition $1<p<p_0(n,\al,\mu)$ includes the one 
$1<p<p_G'(n,\al,\mu)$. 
This means that if $\al$ is close to $1$, then the better region $1<p<p_0(n,\al,\mu)$ than $1<p<p_G'(n,\al,\mu)$ appears for $0\le \mu<1$, as shown by the region (O) in Fig.~\ref{Fig2}. 
Hence, the exponent $p_G'(n,\al,\mu)$ cannot always be the critical exponent for the global existence of solutions.

If we compare the two upper bounds \eqref{subcls-21} and \eqref{subcls-22}, 
then in the region (G) in Fig.~\ref{Fig1} and Fig.~\ref{Fig2}, 
\eqref{subcls-21} is better than \eqref{subcls-22} while in the region (O) in Fig.~\ref{Fig2}, 
this relation becomes reverse.

\subsection{Case $\al \ge 1$}

We next consider the same problem for the case $\al\ge 1$. 

\begin{theo}
Let $n\ge 2, \; 0\le \al<1, \; \mu\ge 0$ and $1<p<1+1/\mu$.
Assume that $u_0\in C^2(\R^n)$ and $u_1\in C^1(\R^n)$ are nontrivial and satisfy $u_1(x)\ge 0$, $\mbox{\rm supp }u_0, \mbox{\rm supp }u_1\subset \{|x|\le R\}$ with $R>0$. 
Suppose that the problem \eqref{Prob0} with \eqref{data0} has  a classical solution $u\in C^2([1,T)\times\R^n)$.  
Then, $T<\infty$ and there exists a constant $\ep_0>0$ depending on 
$p,\al,\mu,R,u_0,u_1$ such that $T_\ep$ has to satisfy 
\begin{align}
&T_\ep^{1-\mu(p-1)}(\ln T_\ep)^{-n(p-1)}\le C\ep^{-(p-1)} &&\mbox{if }\al=1, 
\label{ae1}
\\
&T_\ep\le C\ep^{-(p-1)/\{1-\mu(p-1)\}} &&\mbox{if }\al>1
\label{ag1}
\end{align}
for $0<\ep\le \ep_0$, 
where $C>0$ is a constant independent of $\ep$. 
\label{th23}
\end{theo}
 
\pf 
We first remark that 
the $C^2$-solution $u$ of \eqref{Prob0} and \eqref{data0} has the property of finite speed of propagation, and satisfies
\[
\mbox{supp }u(t,\cdot)\subset \{|x|\le A(t)+R\}, 
\]
where 
\begin{equation}
A(t)=
\begin{cases}
\dsp  \int_1^t s^{-1}ds =\ln t  &\mbox{if }\al=1, \\
 \dsp  \int_1^t s^{-\al}ds =\fa 1{\al-1}(1-t^{1-\al}) &\mbox{if }\al>1, 
\end{cases}  
\label{suppu1}
\end{equation}
provided that supp $u_0$, supp $u_1\subset \{|x|\le R\}$. See \cite{TW3} for its proof. 
\\

Define $F(t)$ by \eqref{int-u_t} and 
let $\al=1$. 
From \eqref{eqF} and \eqref{suppu1}, 
\begin{align}
F'(t)+\fa\mu t F(t) 
 & \ge \fa 1{(A(t)+R)^{n(p-1)}}|F(t)|^p \nonumber \\
 & \ge \fa C{(\ln t)^{n(p-1)}}|F(t)|^p. 
\label{Ho2}
\end{align}
Proceeding as in the proof of Theorem \ref{th21} for the case $0\le \al<1$, 
we have by \eqref{F2},  
\begin{align*}
F(t)&\ge C\ep^p t^{-\mu}\int_1^t s^{\mu-\mu p}(\ln s)^{-n(p-1)}ds  \\
&\ge C\ep^p t^{-\mu(p+1)}(\ln t)^{-n(p-1)}\int_1^t (s-1)^\mu ds.  
\end{align*}
Therefore, we obtain
\begin{equation}
F(t)\ge C\ep^p t^{-\mu(p+1)}(\ln t)^{-n(p-1)}(t-1)^{\mu+1} 
\qquad \mbox{for }t\ge 1. 
\label{ann2}
\end{equation}
Finally, by \eqref{Ho2} and \eqref{ann2}, applying Lemma \ref{lm22} with $q=0$, $a=\mu(p+1)$, $b=r=n(p-1), \; c=\mu+1$ and $A_0=C\ep^p$, we obtain the desired result \eqref{ae1} since
\begin{equation}
M =(p-1)(1 - \mu p)+1=p(1-\mu(p-1))>0. 
\label{Mp2}
\end{equation}

On the other hand, 
if $\al>1$, then  
from \eqref{eqF} and \eqref{suppu1}, 
\begin{equation}
F'(t)+\fa\mu t  F(t) \ge C |F(t)|^p. 
\label{Ho3}
\end{equation}
Proceeding as before, by \eqref{F2},  
we have 
\begin{align*}
F(t)&\ge C\ep^p t^{-\mu}\int_1^t s^{\mu-\mu p}ds  \\
&\ge C\ep^p t^{-\mu(p+1)}\int_1^t (s-1)^\mu ds. 
\end{align*}
Therefore, we obtain
\begin{equation}
F(t)\ge C\ep^p t^{-\mu(p+1)}(t-1)^{\mu+1} 
\qquad \mbox{for }t\ge 1. 
\label{ann3}
\end{equation}
Finally, by \eqref{Ho3} and \eqref{ann3}, applying Lemma \ref{lm22} with $q=0$, $a=\mu(p+1)$, $b=r=0, \; c=\mu+1$ and $A_0=C\ep^p$, we obtain the desired result \eqref{ag1} since \eqref{Mp2}. 
This completes the proof of Theorem \ref{th23}.  
\hfill\qed

\vspace{5mm}

\section{Space derivative nonlinearity.}

\setcounter{equation}{0}

In this section we consider the problem \eqref{Prob0x}. 
Let $u_0$ and $u_1$ be nonnegative and  satisfy supp $u_0$, supp $u_1\subset \{|x|\le R\}$ with $R>0$. 

We prepare several basic inequalities which will be used repeatedly. 
Let 
\[
F(t)=\int u(t,x)dx. 
\]
Then integrating equation \eqref{Prob0x} over $\R^n$ and using Poincar\'e's and H\"older's inequalities imply that 
\begin{align}
F''(t)+\fa\mu t F'(t) &=\int |\nabla_x u|^p dx  
   \label{eqF1} \\
 &\ge  \fa 1{(A(t)+R)^p}\int|u|^pdx  \label{eqF2}\\ 
  &\ge  \fa 1{(A(t)+R)^{p+n(p-1)}}|F(t)|^p. \label{eqF3} 
\end{align}

On the other hand, mutiplying \eqref{eqF1} by $t^\mu$
 and integrating imply
\begin{align}
t^\mu F'(t)-F'(1)&= \int_1^t s^\mu\int |\nabla_x u|^p dxds. \label{impr}
\intertext{Since $F'(1)>0$ by assumption,}  
 F'(t)& \ge t^{-\mu}\int_1^t s^\mu\int |\nabla_x u|^p dxds. \nonumber
\end{align} 
Integrating again, we have from $F(1)>0$ by assumption, 
\begin{equation}
F(t)\ge \int_1^t \tau^{-\mu}\int_1^\tau s^\mu\int |\nabla_x u|^p dxdsd\tau  \qquad 
\mbox{for }t\ge 1. 
\label{Fup}
\end{equation}

\noindent
\subsection{Case $0\le \al<1$}

We call wavelike and heatlike cases if a blow-up condition is concerned with exponents similar to the Strauss and Fujita ones, respectively.

\subsubsection{\large Wavelike and subcritical case}

Let $p_c'(n,\al,\mu)$ be the positive root of the equation
\footnote{
For the equation 
$u_{tt}-\fa 1{t^{2\al}}\Delta u+\fa\mu{t}u_t=|u|^p$,  
the critical exponents as blow-up condtions are 
$p_c(n,\al,\mu)$ and $p_F(n,\al)$, 
where 
$p_c(n,\al,\mu)$ is the positive root of 
\[
\gamma(n,p,\al,\mu)
=-p^2\lp{n-1+\fa{\mu-\al}{1-\al}}+p\lp{n+1+\fa{\mu+3\al}{1-\al}}+2=0, 
\]
and $p_F(n,\al)=1+2/\{n(1-\al)\}$. 
We remark that if $\mu=\al=0$, then $p_c(n,0,0)$ and $p_F(n,0)$ coincide with the Strauss and Fujita exponents, respectively. 
See \cite{TW1,TW2}. 
} 
\begin{align}
\gamma'(n,p,\al,\mu)
&\equiv
-p^2\lp{n+1+\fa{\mu-\al}{1-\al}}+p\lp{n+1+\fa{\mu+3\al}{1-\al}}+2=0
\label{gm'} 
\end{align}
and let 
\[
p_0'(n,\al,\mu)=1+\fa{1+\al}{(n+1)(1-\al)+\mu-1}. 
\]

\begin{theo}
Let $n\ge 2, \; 0\le \al<1, \; \mu\ge 0$ and $1<p<p_c'(n,\al,\mu)$ or 
$1<p<p_0'(n,\al,\mu)$.
Assume that $u_0\in C^2(\R^n)$ and $u_1\in C^1(\R^n)$ are nonnegative, nontrivial and $\mbox{\rm supp }u_0, \mbox{\rm supp }u_1\subset \{|x|\le R\}$ with $R>0$. 
Suppose that the problem \eqref{Prob0x} has  a classical solution $u\in C^2([1,T)\times\R^n)$.  
Then, $T<\infty$ and there exists a constant $\ep_0>0$ depending on 
$p,\al,\mu,R,u_0,u_1$ such that $T_\ep$ has to satisfy 
\begin{align}
T_\ep&\le C\ep^{\fa{-2p(p-1)}{(1-\al)\gamma'(n,p,\al,\mu)}} &&\mbox{ if } 1<p<p_c'(n,\al,\mu), 
\label{subcls-31} \\
T_\ep&\le C\ep^{-\fa{p-1}{\{1-\mu-(n+1)(1-\al)\}(p-1)+1+\al}} &&\mbox{ if }
1<p< p_0'(n,\al,\mu)
\label{subcls-32}
\end{align}
for $0<\ep\le \ep_0$, 
where $C>0$ is a constant independent of $\ep$. 
\label{th31}
\end{theo}
    
\noindent
{\bf Remark}
(1) The upper bound of the lifespan in \eqref{subcls-32} is better than that in \eqref{subcls-31} if 
$1<p< 2(1-\al)/\{(n+1)(1-\al)+\mu+\al-2\}$ since  
\[
\{1-\mu-(n+1)(1-\al)\}(p-1)+1+\al> \fa{1-\al}{2p}\gamma'(n,p,\al,\mu). 
\]
Note that in this case, if $p<p_c'(n,\al,\mu)$, then $\{1-\mu-(n+1)(1-\al)\}(p-1)+1+\al>0$. \\ 
(2) If $p< 2(1-\al)/\{(n+1)(1-\al)+\mu+\al-2\}$, then the condition $\mu<1$ is necessary since $p>1$ and $n\ge 2$. 
Hence, the theorem in this case is not applied to the original equation \eqref{ore} for $|\nb_x u|^p$ in the FLRW spacetime since $\mu=2/(1+w)\ge 1$. This is covered later in Section 4 in more detail. 

\pf 
We have proved in \cite{TW1}, by choosing large $T_1>0$, 
\[
\int |u|^p dx  \ge C\ep^p t^{-(\mu-\al)p/2+(1-\al)(n-1)(1-p/2)}
\]
for $t\ge T_1$. 
By \eqref{suppu0} and \eqref{eqF2}, we obtain
\begin{align}
\int |\nabla_x u|^p dx & \ge C\ep^p t^{-(\mu-\al)p/2+(1-\al)(n-1)(1-p/2)-p(1-\al)} \qquad \mbox{for }t\ge T_1. 
\label{up}
\end{align}
Using \eqref{Fup} implies  
\begin{align}
F(t)&\ge C\ep^p \int_{T_1}^t \tau^{-\mu}\int_{T_1}^\tau s^{\mu-(\mu-\al)p/2+(1-\al)(n-1)(1-p/2)-p(1-\al)}dsd\tau \nonumber \\
&\ge C\ep^p \int_{T_1}^t \tau^{-\mu(1+ p/2)-(1-\al)(n-1)p/2-p(1-\al)}
  \int_{T_1}^\tau (s-T_1)^{\mu+\al p/2+(1-\al)(n-1)}dsd\tau \nonumber\\
&\ge C\ep^p t^{-\mu(1+ p/2)-(1-\al)(n-1)p/2-p(1-\al)}\int_{T_1}^t (\tau-T_1)^{\mu+\al p/2+(1-\al)(n-1)+1}d\tau. \nonumber \\
\intertext{Therefore, }
F(t)&\ge C\ep^p t^{-\mu(1+ p/2)-(1-\al)(n-1)p/2-p(1-\al)}(t-T_1)^{\mu+\al p/2+(1-\al)(n-1)+2} 
\qquad \mbox{for }T\ge T_1. 
\label{an}
\end{align}
We also have
\begin{equation}
F''(t)+\fa\mu t F'(t) \ge \fa C{(t+R)^{(p+n(p-1))(1-\al)}}|F(t)|^p
\label{eqF4}
\end{equation}
by \eqref{suppu0} and \eqref{eqF3}. 
The following lemma is proved in \cite{TW1}. 

\begin{lem} 
Let $p>1, \;a\ge 0, \;b>0, \;q>0, \; \mu\ge 0$ and 
\[
M\equiv (p-1)(b-a)-q+2>0. 
\]
Let $T\ge T_1>T_0\ge 1$. 
Assume that $F\in C^2([T_0,T))$ satisfies the following three conditions:  
\begin{align*}
(i) \quad & F(t) \ge A_0t^{-a}(t-T_1)^b \qquad \mbox{for }t\ge T_1 \\
(ii) \quad &F''(t) +\fa{\mu F'(t)}{t}\ge A_1(t+R)^{-q}|F(t)|^p \quad \mbox{for }t\ge T_0	 \\
(iii) \quad &F(T_0)\ge 0, \quad F'(T_0)>0,  
\end{align*}
where $A_0, A_1$ and $R$ are positive constants. Then 
$T$ has to satisfy
\[
T<C A_0^{-(p-1)/M},  
\]
where $C$ is a constant depending on $R,A_1,\mu,p,q,a$ and $b$. 
\label{lm32}
\end{lem}

From \eqref{an} and \eqref{eqF4}, applying Lemma \ref{lm32} with $q=(p+n(p-1))(1-\al)$, $a=\mu(1+ p/2)+(1-\al)(n-1)p/2+p(1-\al)$, $b=\mu+\al p/2+(1-\al)(n-1)+2$, and $A_0=C\ep^p$, we obtain the desired result since 
\begin{align*}
M&=(p-1)\lb{2-\fa{\mu-\al}2p+(1-\al)\lp{(n-1)\lp{1-\fa p2}-p}}-(p+n(p-1))(1-\al)+2 \\
&=\fa{1-\al}2\gamma'(n,p,\al,\mu)>0.
\end{align*}

It remains to prove \eqref{subcls-32}. 
From \eqref{impr}, we have $F'(t)\ge t^{-\mu}F'(1)$, hence, integrating 
and using $F(1) > 0$ by assumption imply 
\[
F(t)\ge F'(1)t^{1-\mu}=C\ep t^{1-\mu}  \quad \mbox{for }t\ge T'
\]
with some $T'>1$. 
By \eqref{suppu0} and \eqref{eqF3}, 
\begin{equation}
\int |\nabla_x u|^p dx  \ge C\ep^p t^{-(p+n(p-1))(1-\al)+p(1-\mu)} \qquad \mbox{for }t\ge T'. 
\label{up1}
\end{equation}
From \eqref{Fup}, 
\begin{align}
F(t) & \gtrsim \ep^p \int_{T'}^t \tau^{-\mu}\int_{T'}^\tau s^{\mu-(p+n(p-1))(1-\al)+p(1-\mu)}dsd\tau \nonumber \\
&\gtrsim \ep^p \int_{T'}^t \tau^{-(p+n(p-1))(1-\al)-p\mu}\int_{T'}^\tau s^pds d\tau \nonumber \\
&\gtrsim \ep^p \int_{T'}^t \tau^{-(p+n(p-1))(1-\al)-p\mu}(\tau-T')^{p+1} d\tau \nonumber \\
&\ge C\ep^p t^{-(p+n(p-1))(1-\al)-p\mu}(t-T')^{p+2} \qquad \mbox{for }t\ge T'. 
\label{imprFlow}
\end{align}
From \eqref{eqF4} and \eqref{imprFlow}, 
applying Lemma \ref{lm32} with $q=(p+n(p-1))(1-\al)$, $a=(p+n(p-1))(1-\al)+p\mu$, $b=p+2$ and $A_0=C\ep^p$, we obtain the desired result. 
This completes the proof of Theorem \ref{th31}. 
\hfill\qed

In Theorem \ref{th31} the coefficient $n+1+(\mu-\al)/(1-\al)$ of $p^2$ in \eqref{gm'} is positive since it is assumed that the positive root $p_c'(n,\al,\mu)$ exists. 
If $\mu<1$ and $3/4\le\al<1$, then $n+1+(\mu-\al)/(1-\al)\le 0$ can happen; hence,  
$\gamma'(n,p,\al,\mu)>0$ for all $p>1$ and the following corollary holds: 

\begin{coro} 
Let 
$n\ge 2, \; 0\le \al<1, \; \mu\ge 0, \; n+1+(\mu-\al)/(1-\al)\le 0$ and $p>1$. 
Under the assumptions on the initial data of Theorem \ref{th31}, there holds  \eqref{subcls-31}. 
\label{co33}
\end{coro}

\noindent
{\bf Remark} One cannot apply Corollary \ref{co33} to the original equation \eqref{ore} for $|\nb_x u|^p$ in the FLRW spacetime since $n+1+(\mu-\al)/(1-\al)>0$. This is covered later in Section 5 in more detail.

\subsubsection{\large Wavelike and critical case }

We next consider the critical case $p=p_c'(n,\al,\mu)$. 
  
\begin{theo}
 Let $0\le \al<1$ for $n\ge 3$ and $2/7<\al <1$ for $n=2$, and let $\mu\ge 0$ and 
 \begin{align}
 p=p_c'(n,\al,\mu)>
 \begin{cases}
\dsp  p_F'(n,\al)=1+\fa{1+\al}{(n+1)(1-\al)} & \mbox{if }n\ge 3, \\
\max\{p_F'(2,\al),2\} & \mbox{if }n=2. 
\end{cases}
\label{pcf}
\end{align} 
Assume that $u_0\in C^2(\R^n)$ and $u_1\in C^1(\R^n)$ are nonnegative, nontrivial and $\mbox{\rm supp }u_0, \mbox{\rm supp }u_1
\subset \{|x|\le R\}$ with some $0<R\le 1/(2(1-\al))$. 
Suppose that the problem \eqref{Prob0x} has  a classical solution $u\in C^2([1,T)\times\R^n)$.  
Then, $T<\infty$, and there exists a constant $\ep_0>0$ depending on 
$p,\al,\mu,R,u_0,u_1$ such that $T_\ep$ has to satisfy 
\[
T_\ep\le \exp(C\ep^{-p(p-1)}) 
\]
for $0<\ep\le \ep_0$, 
where $C>0$ is a constant independent of $\ep$. 
\label{th34}
\end{theo}

\noindent
{\bf Remark} 
(1) In case $n=2$, we note that 
\[
\max\{p_F'(2,\al),2\}=
\begin{cases}
p_F'(2,\al) & \mbox{if }\al\ge 1/2, \\
2 & \mbox{if }\al < 1/2. 
\end{cases}
\]
(2) When $n=2$, the condition $\al>2/7$ is necessary for the existence of $p$ satisfying \eqref{pcf}. \\
(3) If $p_F'(n,\al)=p_c'(n,\al,\mu)$, then 
\[
\mu=\mu^\ast\equiv (n+1)(1-\al)+\al-
\fa{2(n+1)(1-\al)^2}{(n+1)(1-\al)+1+\al}. 
\]
(4) The case $p_c'(n,\al,\mu)<p\le p_F'(n,\al)$ is considered in the next sub-subsection.  

\pf 
Let
\begin{align}
\lm_\eta(t)&=\lm(\eta t)=(\eta t)^{(1-\mu)/2}K_\nu\lp{\fa 1{1-\al}(\eta t)^{1-\al}}, 
\nonumber \\
\psi_\eta(t,x)&=\lm_\eta(t)\int_{|\omega|=1}e^{\eta^{1-\al}x\cdot \omega}dS_\omega, 
\label{psieta} 
\end{align}
where $K_\nu(t)$ is the modified Bessel function given by \eqref{Bessel}. 
We now define a test function $\phi_q(t,x)$ by 
\begin{equation}
\phi_q(t,x)=\int_0^1 \psi_\eta(t,x)\eta^{q-1+\mu}d\eta.
\label{pqdef}
\end{equation}
Let $q$ satisfy 
\begin{align}
&q>-\fa{\mu+\al}2  \label{cd1}\\
\intertext{and }
&q+\fa{\mu-1}2-(1-\al)|\nu|>-1.  \label{cd2}
\end{align}
It is proved in \cite{TW1} that the function $\phi_q(t,x)$ satisfies the following properties:

\begin{lem}
Let $\phi_q(t,x)$ be defined by \eqref{pqdef}. Assume that 
$q$ satisfies \eqref{cd1} and \eqref{cd2}. 
\begin{enumerate}[(i)]
\item
Then,  
there exists a $T_0>0$ such that 
$\phi_q$ satisfies
\[
\phi_q(t,x) \sim 
\begin{cases}
t^{(-\mu+\al)/2}(t^{1-\al}+|x|)^{-\lp{q+\fa{\mu+\al}2}/(1-\al)} \\
\hspace{6cm} \lp{-\fa{\mu+\al}2<q<\fa{(n-1)(1-\al)-(\mu+\al)}2}, \\
  & \\
t^{(-\mu+\al)/2} (t^{1-\al}+|x|)^{-(n-1)/2} 
(t^{1-\al}-(1-\al)|x|)^{(n-1)/2-\lp{q+\fa{\mu+\al}2}/(1-\al)} \\
 \hspace{7cm}\lp{q>\fa{(n-1)(1-\al)-(\mu+\al)}2},
\end{cases}
\]
for $t\ge T_0$ and $|x|\le(t^{1-\al}-1)/(1-\al)+R$ with $0<R\le 1/(2(1-\al))$. 
\item  Moreover, if $q+1-\al>\{(n-1)(1-\al)-(\mu+\al)\}/2$, 
then there exists a $T_1>0$ such that 
$\phi_q$ satisfies
\[
\p_t\phi_q(t,x) \sim  t^{-(\mu+\al)/2}(t^{1-\al}+|x|)^{-(n-1)/2} 
(t^{1-\al}-(1-\al)|x|)^{(n-1)/2-\lp{q+\fa{\mu+\al}2}/(1-\al)-1} 
\]
for $t\ge T_1$ and $|x|\le(t^{1-\al}-1)/(1-\al)+R$ with $0<R\le 1/(2(1-\al))$. 
\end{enumerate}
\label{lm35} 
\end{lem}

We now prove the following key lemma to prove the theorem. 

\begin{lem}
Assume that $u,u_0$ and $u_1$ satisfy the conditions in Theorem \ref{th34}. 
Let $n\ge 2, \;  \nu=(\mu-1)/(2(1-\al))$ and 
\[
q=\fa{(n-1)(1-\al)-(\mu+\al)}2 -\fa{1-\al}p, 
\]
and let $p$ satisfy \eqref{pcf}. 
Define 
\[
G(t)=\int_1^t(t-\tau)\tau^{1+\mu}\int |\nabla_x u|^p\phi_q (\tau,x)dxd\tau. 
\]
Then, $G(t)$ satisfies 
\[
G'(t)\ge C(\ln t)^{1-p}t\lp{\int_1^t \tau^{-3}G(\tau)d\tau}^p
\qquad \mbox{for }t\ge T_2 
\]
with some $T_2$ sufficiently large, 
where $C$ is a constant independent of $\ep$. 
\label{lm36}
\end{lem}

\pf 
We first verify that $q$ satisfies the required conditions \eqref{cd1} and \eqref{cd2} to use Lemma \ref{lm35}. 
We claim that there holds 
\begin{equation}
-\fa{\mu+\al}2<q<\fa{(n-1)(1-\al)-(\mu+\al)}2. 
\label{q1}
\end{equation}
The second inequality is clearly true. 
To show the first inequality or \eqref{cd1}, we remark that \eqref{pcf} implies after some calculation $p>2/(n-1)$, which is equivalent to $q>-(\mu+\al)/2$. 

We can also show that $q$ satisfies \eqref{cd2} with $\nu=(\mu-1)/(2(1-\al))$, {\it i.e., }$q>-\min\{\mu,1\}$. 
In fact, 
(i) the assumption $p>p_F'(n,\al)$ is equivalent to $q>-1$ since 
\begin{equation}
q=n(1-\al)-\fa 2{p-1}-1+p'(1-\al), \quad \fa 1p+\fa 1{p'}=1  \qquad \mbox{if }p=p_c', 
\label{q-eq}
\end{equation}
(ii) the critical case $p=p_c'$ satisfies $\gamma'(n,p,\al,\mu)=0$ in \eqref{gm'} and this  equality yields 
\begin{align*}
p&=\fa{(n+1)(1-\al)+\mu+3\al}{(n+1)(1-\al)+\mu-\al}+\fa{2(1-\al)}{p\lb{(n+1)(1-\al)+\mu-\al}} \\
&>\fa{(n+1)(1-\al)+\mu+\al}{(n+1)(1-\al)+\mu-\al} \\
&=1+\fa{1+\al}{(n+1)(1-\al)+\mu-1}, 
\end{align*}
which is equivalent to $q>-\mu$ by \eqref{q-eq}.  

In addition, $q$ satisfies the condition of Lemma \ref{lm35} (ii) since  
\begin{align}
q&<\fa{(n-1)(1-\al)-(\mu+\al)}2 \nonumber\\
&<\fa{(n-1)(1-\al)-(\mu+\al)}2+(1-\al)\lp{1-\fa 1p}
=q+1-\al. 
 \label{q2}
\end{align}   

Let us now show the inequality of the lemma. 
Mutiplying eq. in \eqref{Prob0x} by a test funtion $\phi(t,x)$ and $t^\mu$, and 
integrating over $\R^n$, we have 
\begin{equation}
\fa d{dt}\int t^\mu(u\phi)_tdx-2\fa d{dt}\int t^\mu u\phi_t dx
+\int t^\mu u\lp{\phi_{tt}-\fa 1{t^{2\al}}\Delta\phi+\fa\mu t\phi_t}dx
=\int t^\mu |\nabla_x u|^p\phi dx. 
\label{3-9}
\end{equation}
As shown in \cite{TW1}, the test function $\phi_q$ given in \eqref{pqdef} satisfies
\begin{equation}
\lp{\p_t^2-\fa 1{t^{2\al}}\Delta +\fa\mu t\p_t}\phi_q(t,x)=0. 
\label{phq0} 
\end{equation}
Applying \eqref{3-9} with $\phi=\phi_q$ and \eqref{phq0}, we have 
\[
\fa{d^2}{dt^2}\int t^\mu u\phi_q dx-\mu\fa d{dt}\int  t^{\mu-1} u\phi_q dx
-2\fa d{dt}\int t^\mu u\p_t\phi_q dx 
=\int t^\mu |\nabla_x u|^p\phi_q dx. 
\]
Moreover, integrating over $[1,t]$ three times, we obtain 
\begin{align*}
&\int_1^t \tau^\mu\int u\phi_qdxd\tau -\mu\int_1^t(t-\tau)\tau^{\mu-1}\int u\phi_q dxd\tau -2\int_1^t(t-\tau)\tau^\mu\int u\p_\tau\phi_q dxd\tau  \nonumber\\
=&C_{data}(t)
+\fa 12\int_1^t (t-\tau)^2\tau^\mu\int  |\nabla_xu|^p\phi_q dxd\tau, 
\end{align*}
where 
\[
C_{data}(t)
=\ep (t-1)\int u_0(x)\phi_q(1,x)dx
+\fa \ep 2(t-1)^2 \int (u_1(x)\phi_q(1,x)-u_0(x)\p_t\phi_q(1,x))dx. 
\]
We note that $\p_t\phi_q(1,x)\le 0$, which is shown in \cite{TW1}. Hence, by positivity assumption on $u_0$ and $u_1$, it holds that $C_{data}(t)\ge 0$ for $t\ge 1$. Thus, 
\begin{align}
&\int_1^t \tau^\mu\int u\phi_q(\tau,x)dxd\tau -\mu\int_1^t(t-\tau)\tau^{\mu-1}\int u\phi_q(\tau,x) dxd\tau \nonumber \\
& \hspace{7cm}-2\int_1^t(t-\tau)\tau^\mu\int u\p_\tau\phi_q(\tau,x) dxd\tau  \nonumber\\
&\ge \fa 12\int_1^t (t-\tau)^2\tau^\mu\int  |\nabla_xu|^p\phi_q(\tau,x) dxd\tau. 
\label{F1}
\end{align}
Since
\[
G'(t)=\int_1^t \tau^{1+\mu}\int |\nabla_xu|^p\phi_q dxd\tau \quad \mbox{and } \quad  
G''(t)= t^{1+\mu}\int |\nabla_xu|^p\phi_q dxd\tau,  
\]
the right-hand side of \eqref{F1} becomes
\begin{align}
\fa 12\int_1^t (t-\tau)^2\tau^\mu\int  |\nabla_xu|^p\phi_q dxd\tau
=\fa 12\int_1^t (t-\tau)^2\tau^{-1} G''(\tau)d\tau 
=\int_1^t t^2\tau^{-3}G(\tau)d\tau. 
\label{esRHS}
\end{align} 

We now estimate the left-hand side of \eqref{F1}. 
By Poincar\'e's and H\"older's inequalities, the first integral is estimated by
\begin{align*}
I&\equiv \int_1^t \tau^\mu\int u\phi_q(\tau,x)dxd\tau \\  
 &\le \int_1^t \tau^\mu \|\phi_q(\tau)\|_{L^\infty(|x|\le A(\tau)+R)}
    \int_{|x|\le A(\tau)+R} |u(\tau,x)|dxd\tau \\   
  &\le \int_1^t \tau^\mu \|\phi_q(\tau)\|_{L^\infty(|x|\le A(\tau)+R)}(A(\tau)+R)
    \int_{|x|\le A(\tau)+R} |\nabla_x u(\tau,x)|dxd\tau \\    
  &\le \int_1^t \tau^\mu \|\phi_q(\tau)\|_{L^\infty(|x|\le A(\tau)+R)}(A(\tau)+R)\\
  &\qquad \cdot \lp{\int_{|x|\le A(\tau)+R} |\nabla_x u(\tau,x)|^p\phi_q(\tau,x)dx}^{1/p}
    \lp{\int_{|x|\le A(\tau)+R}\phi_q(\tau,x)^{1-p'}dx}^{1/p'}d\tau \\   
&\le\lp{\int_1^t \tau^{1+\mu}\int |\nabla_xu|^p\phi_q(\tau,x)dxd\tau}^{1/p}\\
&\qquad \cdot\lp{\int_1^t \tau^{1+\mu-p'}\|\phi_q(\tau)\|_{L^\infty(|x|\le A(\tau)+R)}^{p'}(A(\tau)+R)^{p'}\int_{|x|\le A(\tau)+R}\phi_q(\tau,x)^{1-p'}dxd\tau }^{1/p'}.  
\end{align*}
We remark that $q$ satisfies \eqref{q1}. 
Applying Lemma \ref{lm35} to the last integral above, we have
\begin{align*}
&\lp{\int_1^{T_0}+\int_{T_0}^t} \tau^{1+\mu-p'}\|\phi_q(\tau)\|_{L^\infty(|x|\le A(\tau)+R)}^{p'}(A(\tau)+R)^{p'}\int_{|x|\le A(\tau)+R}\phi_q(\tau,x)^{1-p'}dxd\tau  \\
 \lsm& C_{T_0}+\int_{T_0}^t \tau^{1+\mu-p'+p'(-\mu+\al)/2-p'(q+(\mu+\al)/2)+p'(1-\al)} \\ 
  & \qquad\qquad  \int_{|x|\le A(\tau)+R} \lb{\tau^{(-\mu+\al)/2}(\tau^{1-\al}+|x|)^{-(q+(\mu+\al)/2)/(1-\al)}}^{1-p'} dxd\tau\\
 \lsm &C_{T_0}+\int_1^t \tau^{n(1-\al)-q-p'/p+p'(1-\al)}d\tau \\
 \lsm & C_{T_0} t^{n(1-\al)-q-p'/p+p'(1-\al)+1},  
\end{align*}
where we have used $A(\tau)+R=(\tau^{1-\al}-1)/(1-\al)+R$, with $R\le 1/(2(1-\al))$. 
Note that if $p=p_c'$, then $q$ satisfies by \eqref{q-eq}
\begin{equation}
n(1-\al)-q-\fa{p'}p+p'(1-\al)=p'. 
\label{cdnqp}
\end{equation}
Thus, we obtain
\begin{equation}
I\lsm G'(t)^{1/p}t^{1+1/p'} \qquad \mbox{for }t\ge T_0. 
\label{esI}
\end{equation}

The second integral on the left-hand side of \eqref{F1} can be estimated as before by 
\begin{align*}
II&\equiv -\mu\int_1^t (t-\tau) \tau^{\mu-1}\int u\phi_q(\tau,x)dxd\tau \\
\lsm &\lp{\int_1^t \tau^{1+\mu}\int |\nabla_xu|^p\phi_q(\tau,x)dxd\tau}^{1/p} \\
  &  \cdot
\lp{\int_1^t (t-\tau)^{p'} \tau^{1+\mu-2p'}\|\phi_q(\tau)\|_{L^\infty(|x|\le A(\tau)+R)}^{p'}(A(\tau)+R)^{p'}\int_{|x|\le A(\tau)+R} \phi_q(\tau,x)^{1-p'}dxd\tau }^{1/p'}. 
\end{align*}
Using Lemma \ref{lm35},
\begin{align*}
II \lsm & G'(t)^{1/p}\lp{C_{T_0}+\int_{T_0}^t (t-\tau)^{p'} \tau^{1-2p'-q+n(1-\al)+p'(1-\al)}d\tau}^{1/p'}. 
\end{align*}
Since $q$ satisfies $n(1-\al)-q+1-2p'+p'(1-\al)=0$ by \eqref{cdnqp}, 
we obtain
\begin{equation}
II\lsm G'(t)^{1/p}t^{1+1/p'} \qquad \mbox{for }t\ge T_0. 
\label{esII}
\end{equation}

We finally estimate the third integral on the left-hand side of \eqref{F1},
\[
III\equiv -2\int_1^t (t-\tau) \tau^\mu\int u\p_\tau\phi_q(\tau,x)dxd\tau. 
\]
Set $T_2=\max\{T_0,T_1\}$ to apply Lemma \ref{lm35} (i) and (ii). 
Proceeding in a similar way as before, we have 
\begin{align*}
III\lsm &\lp{\int_1^t \tau^{1+\mu}\int |\nabla_xu|^p\phi_q(\tau,x)dxd\tau}^{1/p}\\
&\quad \qquad \cdot
\lp{C_{T_2}+\int_{T_2}^t (t-\tau)^{p'} \tau^{1+\mu-p'+p'(1-\al)}\int_{|x|\le A(\tau)+R} \phi_q\lp{\fa{\p_\tau\phi_q}{\phi_q}}^{p'}(\tau,x)dxd\tau }^{1/p'}. 
\end{align*} 
We remark here that $q$ satisfies \eqref{q2}. 
By Lemma \ref{lm35},   
\begin{align*}
&\phi_q\lp{\fa{\p_\tau\phi_q}{\phi_q}}^{p'}
\sim \tau^{(-\mu+\al)/2-\al p'}(\tau^{1-\al}+|x|)^{-(n-1)/2-(p'-1)/p} 
(\tau^{1-\al}-(1-\al)|x|)^{-1}  \quad \mbox{for }t\ge T_2. 
\end{align*}
Hence, 
\begin{align*}
&\tau^{1+\mu-p'+p'(1-\al)}\int_{|x|\le A(\tau)+R} \phi_q\lp{\fa{\p_\tau\phi_q}{\phi_q}}^{p'}(\tau,x)dx \\
\le & \tau^{1+\mu-p'+p'(1-\al)+(-\mu+\al)/2-\al p'-(n-1)(1-\al)/2-(p'-1)(1-\al)/p}\int_{|x|\le A(\tau)+R} 
(\tau^{1-\al}-(1-\al)|x|)^{-1} dx \\
\lsm & \tau^{-(n-1)(1-\al)}\int_0^{A(\tau)+R}(\tau^{1-\al}-(1-\al)r)^{-1}r^{n-1}dr \\
\lsm& \ln\tau,  
\end{align*}
where we note that since $(n-1)(1-\al)/2=q+(\mu+\al)/2+(1-\al)/p$ and \eqref{cdnqp},  
\begin{align*}
&1+\mu-p'+p'(1-\al)+\fa{-\mu+\al}2-\al p'-\fa{(n-1)(1-\al)}2-\fa{(p'-1)(1-\al)}p \\
=&1-p'+p'(1-\al)-\al p'-q-\fa{p'(1-\al)}p \\
=&-(n-1)(1-\al).
\end{align*}
Thus, we obtain 
\begin{align}
III&\lsm G'(t)^{1/p}
\lp{C_{T_2}+\int_{T_2}^t (t-\tau)^{p'} \ln\tau d\tau }^{1/p'} \nonumber \\
  &\lsm  G'(t)^{1/p}t^{1+1/p'}(\ln t)^{1/p'} \qquad \mbox{for }t\ge T_2. 
\label{esIII}  
\end{align}
Combining \eqref{F1}, \eqref{esRHS}, and \eqref{esI}-\eqref{esIII} all together, we obtain the desired inequality. This completes the proof of Lemma \ref{lm36}. 
\hfill\qed

Then the rest of the proof is the same as that of Theorem 2.3 in \cite{TW1}. 
\hfill\qed

\subsubsection{\large Heatlike case}

\begin{theo} 
Let $n\ge 2, \; 0\le \al<1, \; \mu\ge 0$ and $1<p\le p_F'(n,\al)=1+(1+\al)/\{(n+1)(1-\al)\}$.
Assume that $u_0\in C^2(\R^n)$ and $u_1\in C^1(\R^n)$ are nonnegative, nontrivial and $\mbox{\rm supp }u_0, \mbox{\rm supp }u_1\subset \{|x|\le R\}$ with $R>0$. 
Suppose that the problem \eqref{Prob0x} has  a classical solution $u\in C^2([1,T)\times\R^n)$.  
Then, $T<\infty$ and there exists a constant $\ep_0>0$ depending on 
$p,\al,\mu,R,u_0,u_1$ such that $T_\ep$ has to satisfy 
\begin{align}
T_\ep&\le C\ep^{-\fa{p-1}{2-\{n(p-1)+p\}(1-\al)}} && \mbox{if }p<p_F'(n,\al),
\label{subcls-33}\\
T_\ep&\le \exp\lp{C\ep^{-p(p-1)/(p+1)}} && \mbox{if } p=p_F'(n,\al) \mbox{ and }0\le \mu \le 1,
\label{subcls-34}\\
T_\ep&\le \exp\lp{C\ep^{-(p-1)}} && \mbox{if }p=p_F'(n,\al) \mbox{ and }\mu > 1  \nonumber
\end{align}
for $0<\ep\le \ep_0$, 
where $C>0$ is a constant independent of $\ep$.
\label{th37}
\end{theo}

\noindent
{\bf Remark}\quad In the critical case $p=p_F'(n,\al)$, if $0\le \mu \le 1$, then the estimate above is better than that for the case $\mu >1$. 
However, \eqref{subcls-31} and \eqref{subcls-32} are applicable for $p=p_F'(n,\al)$ and 
$0\le \mu \le 1$. These estimates are much better than $T_\ep\le \exp\lp{C\ep^{-p(p-1)/(p+1)}}$ above.   
See Fig.s~\ref{Fig3}, \ref{Fig4} and \ref{Fig5} below. 

\pf Proceeding as in \cite{TW2}, we have 
\begin{equation}
F(t)\ge F(1)=C\ep>0 \qquad \mbox{for }t\ge 1. 
\label{F(1)}
\end{equation}
By \eqref{suppu0} and \eqref{eqF3}, we have
\begin{equation}
\int |\nabla_x u|^p dx  \ge C\ep^p t^{-(p+n(p-1))(1-\al)} \qquad \mbox{for }t\ge 1. 
\label{uph}
\end{equation}
Let $p<p_F'(n,\al)$. 
From \eqref{Fup}, 
\begin{align}
F(t) & \gtrsim \ep^p \int_1^t \tau^{-\mu}\int_1^\tau s^{\mu-(p+n(p-1))(1-\al)}dsd\tau \nonumber \\
&\gtrsim \ep^p \int_1^t \tau^{-\mu-(p+n(p-1))(1-\al)}\int_1^\tau (s-1)^\mu ds d\tau \nonumber \\
&\gtrsim \ep^p \int_1^t \tau^{-\mu-(p+n(p-1))(1-\al)}(\tau-1)^{\mu+1} d\tau \nonumber \\
&\ge C\ep^p t^{-\mu-(p+n(p-1))(1-\al)}(t-1)^{\mu+2} \qquad \mbox{for }t\ge 1. 
\label{imprFlow-h}
\end{align}
From \eqref{eqF4} and \eqref{imprFlow-h}, 
applying Lemma \ref{lm32} with $q=(p+n(p-1))(1-\al)$, $a=\mu+(p+n(p-1))(1-\al)$, $b=\mu+2$ and $A_0=C\ep^p$, we obtain the desired results since
\begin{align*}
M &=(p-1)\lb{2-(p+n(p-1))(1-\al)}-(p+n(p-1))(1-\al)+2\\
&=p\lb{2-(p+n(p-1))(1-\al)}>0. 
\end{align*}
  
Let next $p=p'_F(n,\al)$. Since $(p+n(p-1))(1-\al)=2$, from \eqref{uph}, 
\[
\int |\nabla_x u|^p dx  \ge C\ep t^{-2} \qquad \mbox{for }t\ge 1. 
\]
Hence, by \eqref{Fup}, 
\begin{align}
F(t) & \gtrsim \ep^p \int_1^t \tau^{-\mu}\int_1^\tau s^{\mu-2}dsd\tau \nonumber \\
&\gtrsim \ep^p \int_1^t \tau^{-\mu-2}\int_1^\tau (s-1)^\mu ds d\tau \nonumber \\
&\gtrsim \ep^p \int_2^t \tau^{-1} d\tau \nonumber \\
&\ge C\ep^p \ln\fa t2 \qquad \mbox{for }t\ge 2. 
\label{imprFlow-hc}
\end{align}
The following lemma is proved in \cite{TW2}. 

\begin{lem}
Let $p>1, \;b>0, \; \mu\ge 0$ and $T\ge T_1>T_0\ge 1$. 
Assume that $F\in C^2([T_0,T))$ satisfies the following three conditions:  
\begin{align*}
(i) \quad & F(t) \ge A_0\lp{\ln \fa t{T_1}}^b\qquad \mbox{for }t\ge T_1, \\
(ii) \quad &F''(t) +\fa{\mu F'(t)}{t}\ge A_1(t+R)^{-2}|F(t)|^p \quad \mbox{for }t\ge T_0,	 \\
(iii) \quad &F(T_0)\ge 0, \quad F'(T_0)>0,  
\end{align*}
where $A_0, A_1$ and $R$ are positive constants. Then,  
$T$ has to satisfy
\[
T<
\begin{cases}
\exp\lp{CA_0^{-(p-1)/\{b(p-1)+2\}}} & \mbox{if }\mu\le 1, \\
\exp\lp{CA_0^{-(p-1)/\{b(p-1)+1\}}} & \mbox{if }\mu > 1, 
\end{cases}
\] 
where $C$ is a constant depending on $R,A_1,\mu,p$ and $b$. 
\label{lm38}
\end{lem}

By \eqref{eqF4} with $(p+n(p-1))(1-\al)=2$  and \eqref{imprFlow-hc}, using Lemma \ref{lm38} with $b=1$ and $A_0 = C\ep^p$, we obtain the desired results. 
This completes the proof of Theorem \ref{th37}. 
\hfill\qed

In the end of this subsection, we discuss the blow-up condtions and the estimates of the lifespan in the subcritical cases in Theorems \ref{th31} and \ref{th37}. 

We recall that $p_c'(n,\al,\mu)$ is the positive root of the equation $\gamma'(n,p,\al,\mu)=0$ given in \eqref{gm'}. 
Fig.s~\ref{Fig3}, \ref{Fig4} and \ref{Fig5}
\begin{figure}[h!]
\includegraphics[width=13cm, bb=120 600 500 760, clip]{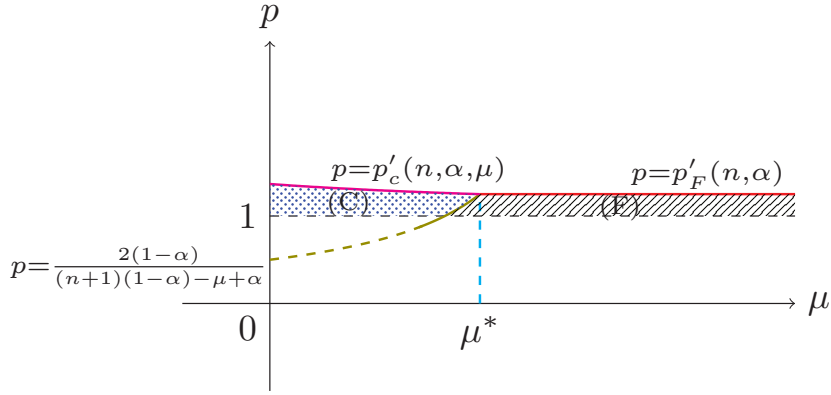}
\caption{\label{Fig3} Range of blow-up conditions in case $n=3$ and $\alpha=0$}
\end{figure}
\begin{figure}[h!]
\includegraphics[width=15cm, bb=120 600 540 770, clip]{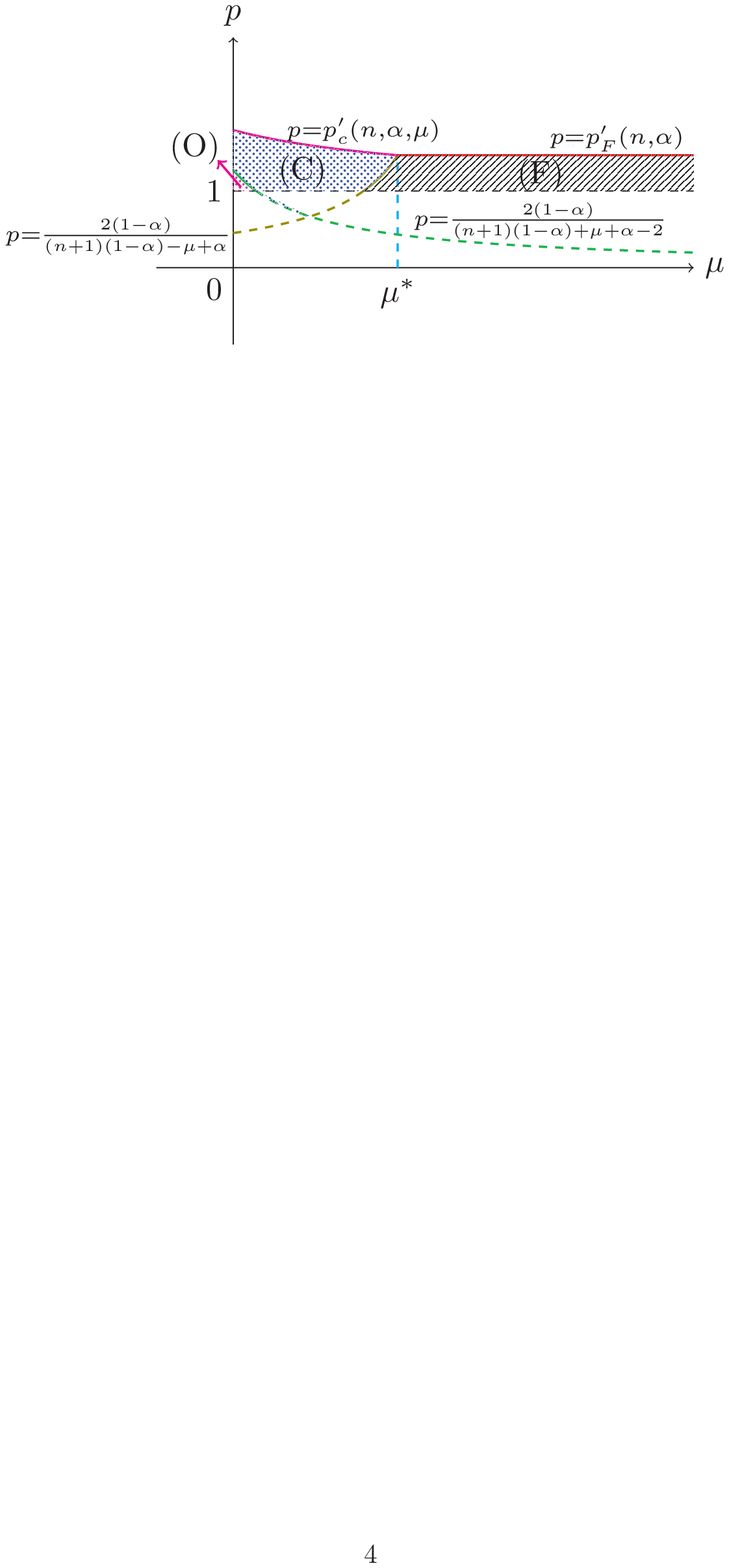}
\caption{\label{Fig4} Range of blow-up conditions in case $n=3$ and $\alpha=0.3$}
\end{figure}
\begin{figure}[h!]
\includegraphics[width=15cm, bb=120 450 580 760, clip]{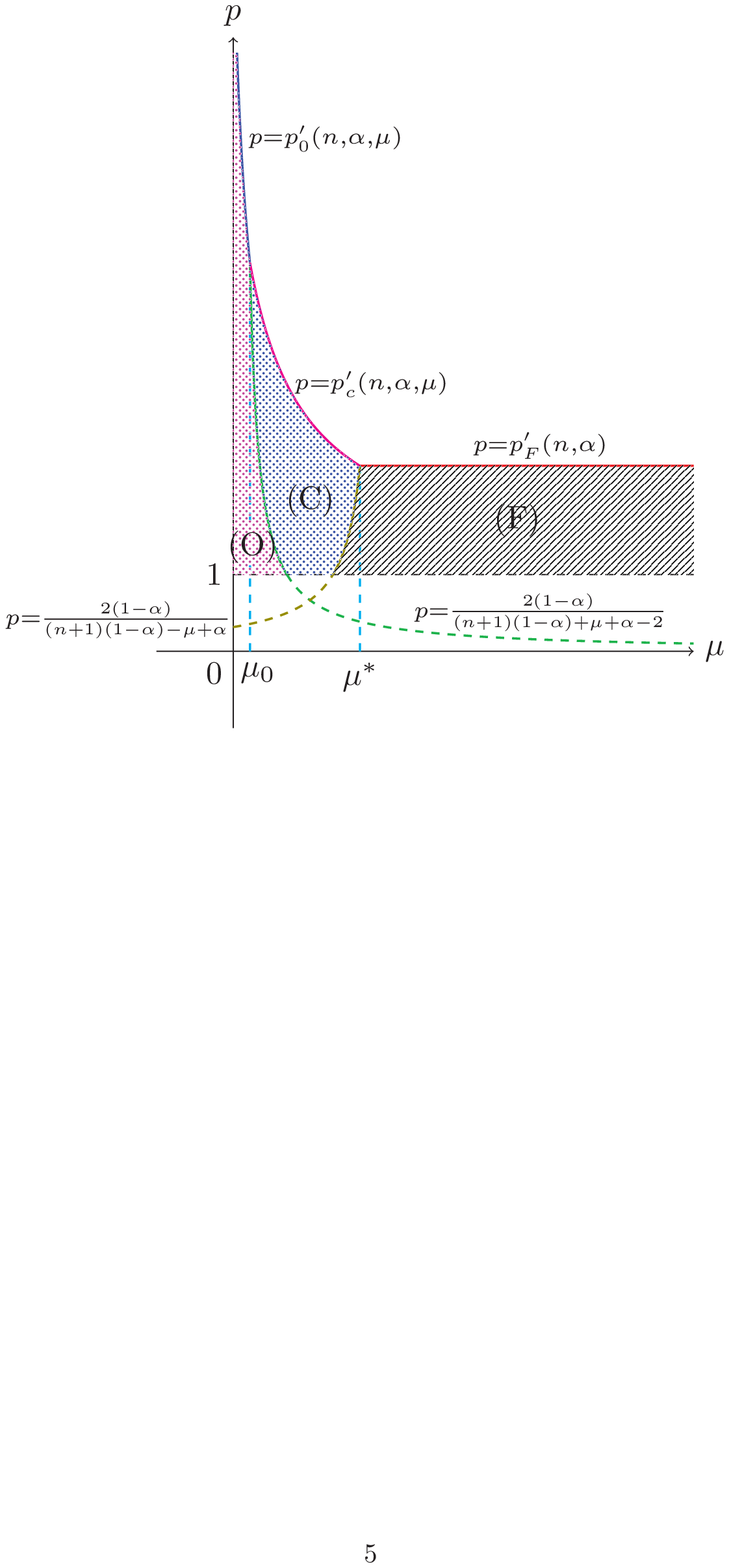}
\caption{\label{Fig5} Range of blow-up conditions in case $n=3$ and $\alpha=0.7$}
\end{figure} 
show the regions of blow-up conditions in the case $\al=0, \; \al=0.3$ and $\al=0.7$, respectively, each for $n=3$.  
We recall that if $p_c'(n,\al,\mu)=p_F'(n,\al)$, then 
\[
\mu=\mu^\ast= (n+1)(1-\al)+\al-\fa{2(n+1)(1-\al)^2}{(n+1)(1-\al)+1+\al}. 
\]
Note that if $p_c'(n,\al,\mu)=p_0'(n,\al,\mu)
$, then 
\[
\mu=\mu_0\equiv -(n-1)(1-\al)+\sr{3\al^2-4\al+2}. 
\]
We easily see that $\mu^\ast >1$ and $\mu_0<1$. 
We also note that $2p(p-1)/\{(1-\al)\gamma'(n,p,\al,\mu)\}=(p-1)/[2-\{n(p-1)+p\}(1-\al)]$ yields $p= 2(1-\al)/\{(n+1)(1-\al)-\mu+\al\}$, and that 
$2p(p-1)/\{(1-\al)\gamma'(n,p,\al,\mu)\}=(p-1)/[\{1-\mu-(n+1)(1-\al)\}(p-1)+1+\al]$ yields $p= 2(1-\al)/\{(n+1)(1-\al)+\mu+\al-2\}$. 

Among the three upper bounds \eqref{subcls-31}, \eqref{subcls-32} and \eqref{subcls-33}, if $1<p< \min\{p_0'(n,\al,\mu), \; 2(1-\al)/\{(n+1)(1-\al)+\mu+\al-2\}\}$, then \eqref{subcls-32} is the best. This is Region (O) shown in Fig.s~\ref{Fig4} and \ref{Fig5}. 

If $\max\{2(1-\al)/\{(n+1)(1-\al)+\mu+\al-2\}, \; 2(1-\al)/\{(n+1)(1-\al)-\mu+\al\}, \; 1\}<p<p_c'(n,\al,\mu)$, then \eqref{subcls-31} is the best (Region (C) in Fig.s~\ref{Fig3}, \ref{Fig4} and \ref{Fig5}). 

On the other hand, 
if $1<p\le 2(1-\al)/\{(n+1)(1-\al)-\mu+\al\}$ and $p<p_F'(n,\al)$, then 
\eqref{subcls-33} is the best (Region (F) in 
Fig.s~\ref{Fig3}, \ref{Fig4} and \ref{Fig5}). 

After some calculation, we see the following facts. 
If $(n-3)/(n-2)\le \al <1$ for $n\ge 3$ and $0\le \al<1$ for $n=2$, 
Region (O) appears in $0\le \mu <1$. 
Moreover, if $\al\ge \max\{0, \; (n^2-2n-1-\sr{n^2-2n-1})/(n^2-2n-2)\}$, 
then $p_0'(n,\al,\mu)>p_c'(n,\al,\mu)>p_F'(n,\al)$ for $0\le \mu< \mu_0$, 
hence the blow-up condition $1<p\le p_0'(n,\al,\mu)$ is the best for $0\le \mu< \mu_0$. 
This is unlike the case of the equation with $|u|^p$-nonlinearity for which 
the blow-up condition is related to only $p_c(n,\al,\mu)$ and $p_F(n,\al)$.

\subsection{Case $\al\ge 1$}

\begin{theo}
 Let $n\ge 2, \; \al\ge 1, \; \mu\ge 0$ and $p>1$.
Assume that $u_0\in C^2(\R^n)$ and $u_1\in C^1(\R^n)$ are nonnegative, nontrivial and $\mbox{\rm supp }u_0, \mbox{\rm supp }u_1\subset \{|x|\le R\}$ with $R>0$. 
Suppose that the problem \eqref{Prob0x} has  a classical solution $u\in C^2([1,T)\times\R^n)$.  
Then, $T<\infty$ and there exists a constant $\ep_0>0$ depending on 
$p,\al,\mu,R,u_0,u_1$ such that $T_\ep$ has to satisfy 
\begin{align*}
&T_\ep^2(\ln T_\ep)^{-(p+n(p-1))}\le C\ep^{-(p-1)} &&\mbox{if }\al=1, \\
&T_\ep\le C\ep^{-(p-1)/2} &&\mbox{if }\al>1
\end{align*}
for $0<\ep\le \ep_0$, 
where $C>0$ is a constant independent of $\ep$. 
\label{th39}
\end{theo}

\pf We first prove the theorem for the case $\al=1$. 
By \eqref{suppu1} and \eqref{eqF3}, 
\begin{equation}
F''(t)+\fa\mu t F'(t) \ge \fa C{(\ln t)^{p+n(p-1)}}|F(t)|^p. \label{eqF5}
\end{equation}
Proceeding as in \cite{TW3}, we have \eqref{F(1)}. 
\[
F(t)\ge F(1)=C\ep>0 \qquad \mbox{for }t\ge 1. 
\]
By \eqref{suppu1} and \eqref{eqF3}, 
\begin{equation}
\int |\nabla_x u|^p dx  \ge C\ep^p (\ln t)^{-(p+n(p-1))} \qquad \mbox{for }t\ge 1. 
\label{uph-2}
\end{equation}
Hence, from \eqref{Fup}, 
\begin{align} 
F(t)&\ge C\ep^p \int_1^t  (\ln \tau)^{-(p+n(p-1))}\tau^{-\mu}(\tau-1)^{\mu+1}d\tau \nonumber \\
 &\ge C\ep^p (\ln t)^{-(p+n(p-1))}t^{-\mu}\int_1^t (\tau-1)^{\mu+1}d\tau \nonumber \\
 &\ge C\ep^p (\ln t)^{-(p+n(p-1))}t^{-\mu}(t-1)^{\mu+2}
 \qquad 
\mbox{for }t\ge 1. 
\label{Fup1} 
\end{align} 
We here use another Kato's lemma. 
We can combine Lemmas 2.3 and 3.3 in \cite{TW3} to obtain the following lemma. 

\begin{lem}
Let $p>1, \;a\ge 0, \;b\ge 0, \; c>0, \;q\ge 0, \; \mu\ge 0$ and 
\[
M\equiv (p-1)(c-a)+2>0. 
\]
Let $T\ge T_1>T_0\ge 1$. 
Assume that $F\in C^1([T_0,T))$ satisfies the following three conditions:  
\begin{align*}
(i) \quad & F(t) \ge A_0t^{-a}(\ln t)^{-b}(t-T_1)^c \qquad \mbox{for }t\ge T_1, \\
(ii) \quad &F''(t) +\fa{\mu}{t}F'(t)\ge A_1(\ln t)^{-q}|F(t)|^p \quad \mbox{for }t\ge T_0,	 \\
(iii) \quad &F(T_0)> 0,   \quad F'(T_0)>0,  
\end{align*}
where $A_0, A_1$ and $R$ are positive constants. Then,  
$T$ has to satisfy
\[
T^{M/(p-1)}(\ln t)^{-b-q/(p-1)}<CA_0^{-1},  
\]
where $C$ is a constant depending on $R,A_1,\mu,p,q,a,b$ and $c$. 
\label{lm310}
\end{lem}

By \eqref{eqF5} and \eqref{Fup1}, applying Lemma \ref{lm310} with $a=\mu$, $b=q=p+n(p-1)$, $c=\mu+2$ and $A_0=C\ep^p$, we obtain the desired result for $\al=1$ since
\[
M=2(p-1)+2=2p>0. 
\]

It remains to prove the theorem for $\al>1$. 
By \eqref{suppu1} and \eqref{eqF3}, 
\begin{align}
F''(t)+\fa\mu t F'(t) \ge C|F(t)|^p. 
\label{3Fp}
\end{align}
As above, we also have 
\[
F(t)\ge C\ep^p t^{-\mu}(t-1)^{\mu+2}
 \qquad 
\mbox{for }t\ge 1. 
\]
Using Lemma \ref{lm310} again, we obtain the desired results. 
This completes the proof of Theorem \ref{th39}. 
\hfill\qed


\section{Wave Equations in FLRW}
\setcounter{equation}{0}

We now apply Theorems \ref{th21}, \ref{th23}, \ref{th31}, \ref{th34}, \ref{th37} and \ref{th39} to the original equation \eqref{ore}
, which is equivalent to \eqref{Prob0} and \eqref{Prob0x} with $\al=2/(n(1+w))$ and $\mu=2/(1+w)$. 
We treat the case $-1<w\le 1$ and $n\ge 2$ so that $\al\ge 1/n$ and $\mu\ge 1$ in \eqref{Prob0} and \eqref{Prob0x}. 
Observe that the cases $-1<w\le 2/n-1$ and $2/n-1<w\le 1$ correspond to accelerating and decelerating expanding universes, respectively.


Consider the equation with the time derivative nonlinear term $|u_t|^p$. 
Denote $p_G'(n,2/(n(1+w)),2/(1+w))$ by $p_G'(n,w)$. 
Fig.~\ref{Fig6} below shows the range of blow-up conditions in terms of $w$ and $p$ in the case $n=3$. 

For $2/n-1<w\le 1$ and $n\ge 2$, applying Theorem \ref{th21} to \eqref{ore} for $|u_t|^p$, we obtain the following upper bounds of the lifespan: 
\[
\begin{cases}
T_\ep\le C\ep^{\fa{-(p-1)}{1-\{n-1+4/(n(1+w))\}(p-1)/2}}  &  \mbox{if } 1<p<p_G'(n,w),  \\
T_\ep\le \exp(C\ep^{-(p-1)}) &  \mbox{if } p=p_G'(n,w),  
\end{cases}
\]
where $C>0$ is a constant independent of $\ep$. 
We note that the estimate \eqref{subcls-22} in Theorem \ref{th21} is not applied to \eqref{ore} since $\mu\ge 1$ in our case. 
See Region (G) in 
Fig.~\ref{Fig6}. 
\begin{figure}[h!]
\includegraphics[width=13cm, bb=100 570 500 760, clip]{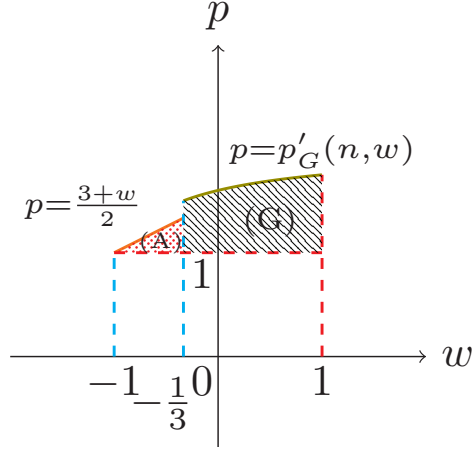}
\caption{\label{Fig6} Range of blow-up conditions in case $n=3$}
\end{figure} 

For $-1<w\le 2/n-1$ and $n\ge 2$, we obtain from Theorem \ref{th23}
\[
\begin{cases}
T_\ep^{1-2(p-1)/(1+w)}(\ln T_\ep)^{-n(p-1)}\le C\ep^{-(p-1)} &  \mbox{if } 1<p<1+\fa 1n \mbox{ and }w=\fa 2n-1,  \\
T_\ep\le C\ep^{-(p-1)/\{1-2(p-1)/(1+w)\}} &  \mbox{if } 1<p<\fa{3+w}2 \mbox{ and }-1<w < \fa 2n-1.   
\end{cases}
\]
See Region (A) in Fig.~\ref{Fig6}. 

From these results, we see that the blow-up range of $p$ in the flat FLRW spacetime is smaller than that in the Minkowski spacetime 
because $p_G'(n,w)<p_G(n)=1+2/(n-1)$.  
Moreover, in the subcritical case $p<p_G'(n,w)$, the lifespan of the blow-up solutions in the FLRW spacetime is longer than that in the Minkowski spacetime since
$\ep^{-(p-1)/\{1-(n-1)(p-1)/2\}}<\ep^{-(p-1)/\{1-(n-1+4/(n(1+w)))(p-1)/2\}}$ for sufficiently small $\ep$. 
Although the critical value has not been established, 
we can say at the present time that global solutions in the FLRW spacetime exist 
more easily than in the Minkowski spacetime.


We next consider the equation with the space derivative nonlinear term $|\nb_x u|^p$. 
We define here $\gamma_0'(n,p,w)$ corresponding to $\gamma'(n,p,\al,\mu)$ in \eqref{gm'}
by 
\[
\gamma_0'(n,p,w)=\lp{1-\fa 2{n(1+w)}}\gamma'\lp{n,p,\fa 2{n(1+w)},\fa 2{1+w}}. 
\]
Then, we obtain
\[
\gamma_0'(n,p,w)=-\lp{n+1-\fa 4{n(1+w)}}p^2+\lp{n+1+\fa 4{n(1+w)}}p+2-\fa 4{n(1+w)}. 
\]
Let $p_c'(n,w)$ be the positive root of the equation $\gamma_0'(n,p,w)=0$.  
We also denote $p_F'(n,2/(n(1+w)))$ by $p_F'(n,w)$. 

For $2/n-1<w\le 1$ and $n\ge 2$, say a decelerated expanding universe, applying Theorems \ref{th31}, \ref{th34} and \ref{th37} to \eqref{ore} for $|\nb_x u|^p$, we obtain the following upper bounds of the lifespan: 
\begin{align}
&T_\ep\le C\ep^{\fa{-2p(p-1)}{\gamma_0'(n,p,w)}}  &&  \mbox{if } 1<p<p_c'(n,w), 
   \label{ls-cw}  \\
&T_\ep\le \exp(C\ep^{-p(p-1)}) &&  \mbox{if } p=p_c'(n,w)>p_F'(n,w)   \nonumber\\
&T_\ep\le C\ep^{\fa{-(p-1)}{2-\{(n+1)p-n\}\{1-2/(n(1+w))\}}}  &&  \mbox{if } 1<p<p_F'(n,w),      \label{ls-Fw} \\
&T_\ep\le \exp(C\ep^{-(p-1)}) &&  \mbox{if } p=p_F'(n,w).  \nonumber  
\end{align}
We note that the estimates \eqref{subcls-32} in Theorem \ref{th31} and \eqref{subcls-34} in Theorem \ref{th37} are not applied to \eqref{ore} since $\mu\ge 1$ in our case, and also that if $n=2$, then $p_F'(2,w)\ge 2$ and $\al=2/(n(1+w))>2/7$ in Theorem \ref{th34} since $w\le 1$. 

For $-1<w\le 2/n-1$ and $n\ge 2$, say an accelerated expanding universe, we obtain from Theorem \ref{th39}
\begin{align}
&T_\ep^2(\ln T_\ep)^{-(n(p+1)-n)}\le C\ep^{-(p-1)} &&  \mbox{if } p>1\mbox{ and }w=\fa 2n-1,  \nonumber\\
&T_\ep\le C\ep^{-(p-1)/2} &&  \mbox{if } p>1\mbox{ and }-1<w < \fa 2n-1.   
\label{ls-aw}
\end{align}
We see that blow-up in finite time can happen to occur for all $p>1$. This is in contrast to the case of decelerated expansion above. 

\begin{figure}[h!]
\includegraphics[width=13cm, bb=50 520 500 760, clip]{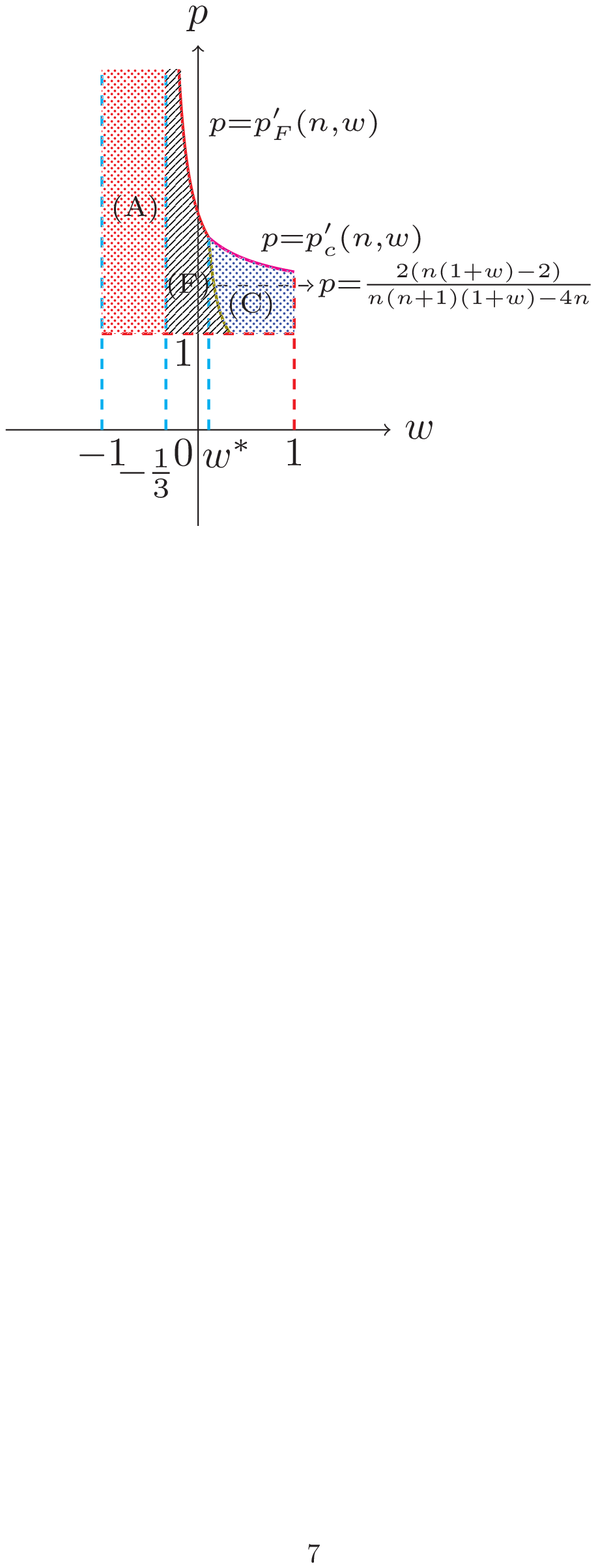}
\caption{\label{Fig7} Range of blow-up conditions in case $n=3$}
\end{figure}

Fig.~\ref{Fig7}
shows the range of blow-up conditions in terms of $w$ and $p$ in the case $n=3$. 
Note that if $p_F'(n,w)=p_c'(n,w)$, then $w$ is the larger root $w^\ast$ of the equation
\[
n^3(n+1)w^2+2n\{n^2(n+1)-(3n+4)(n-1)\}w+n^3(n+1)-2n(3n+4)(n-1)+8(n^2-n-1)=0. 
\]
Region (A) is for the case of the accelerated expanding universe, where the lifespan of blow-up solutions is dominated by \eqref{ls-aw}. 
In contrast, each Region (F) and (C) represents the decelerated expanding one. 
In Region (F) the estimate \eqref{ls-Fw} is better than \eqref{ls-cw}, on the other hand, this relation becomes reverse in Region (C). 

Finally, let us compare the results for the term $|u_t|^p$ with those for $|\nb_xu|^p$, 
especially in the decelerated expanding universe, say (G), (F) and (C). 
We observe that which has a larger blow-up range depends on the value of $w$ for each $n$. 
If $n=2$,  then $\max\{p_F'(2,w),p_c'(2,w)\}>p_G'(2,w)$ for $0=2/n-1<w\le 1$. 
The higher the dimension $n$ becomes, however, the larger the $w$-interval such that $p_G'(n,w)>\max\{p_F'(n,w),p_c'(n,w)\}$ becomes. 

We will treat the remaining case $w=-1$ in future papers.



\end{document}